\documentclass[12pt]{article}
\usepackage{amssymb}
\usepackage{mathrsfs}
\usepackage{amsfonts}
\usepackage{graphicx}
\usepackage{colortbl,dcolumn}
\usepackage{latexsym}
\usepackage{amsbsy,amsthm}
\usepackage{amsmath}
\usepackage{bbm}
\usepackage{subfigure}
\usepackage{psfrag}
\usepackage{float} 
\usepackage{exscale}
\usepackage{relsize}
\usepackage{booktabs}
\usepackage{verbatim}
\usepackage{enumerate}
\usepackage{cite}
\usepackage{amsthm}
\usepackage[top=0.8in, bottom=0.8in, left=0.8in, right=0.8in]{geometry}
\usepackage{hyperref}

\numberwithin{equation}{section}

\newtheorem{assumption}{Assumption}[section]
\newtheorem{remark}[assumption]{Remark}
\newtheorem{lemma}[assumption]{Lemma}
\newtheorem{theorem}[assumption]{Theorem}   
\newtheorem{proposition}[assumption]{Proposition}
\newtheorem{corollary}[assumption]{Corollary}

\hypersetup{colorlinks=true,                         
linkcolor=blue,
anchorcolor=blue,
citecolor=blue}

\allowdisplaybreaks[4] 

\begin{document}
\title
{
    An explicit splitting SAV scheme for the kinetic Langevin dynamics\footnote
    {
    This work is supported by Natural Science Foundation of China (12471394, 12371417, 12071488) and Hunan Basic Science Research Center for Mathematical Analysis (2024JC2002). 
    }
}
\author
{
    Lei Dai\footnotemark[2] \footnotemark[3] , Yingsong Jiang\footnotemark[2] , Xiaojie Wang\footnotemark[2] \footnotemark[4]
}
\maketitle

\footnotetext[2]{School of Mathematics and Statistics, HNP-LAMA, Central South University, Changsha 410083, China.}
\footnotetext[3]{Section of Mathematics, University of Geneva, 1211 Geneva, Switzerland.}
\footnotetext[4]{Corresponding author.
\\
\indent \ Emails: \{dailei, x.j.wang7, yingsong\}@csu.edu.cn; x.j.wang7@gmail.com.}

\begin{abstract}
{
    \rm\small
    The kinetic Langevin dynamics finds diverse applications in various disciplines such as molecular dynamics and Hamiltonian Monte Carlo sampling.
    In this paper, a novel splitting scalar auxiliary variable (SSAV) scheme is proposed for the dynamics, where the gradient of the potential $U$ is possibly non-globally Lipschitz continuous with superlinear growth. As an explicit scheme, the SSAV method is efficient, robust and is able to reproduce the energy structure of the original dynamics. By an energy argument, the SSAV scheme is proved to possess an exponential integrability property, which is crucial to establishing the order-one strong convergence without the global monotonicity condition. Moreover, moments of the numerical approximations are shown to have polynomial growth with respect to the time length. This helps us to obtain weak error estimates of order one, with error constants polynomially (not exponentially) depending on the time length. Despite the obtained polynomial growth, the explicit scheme is shown to be computationally effective for the approximation of the invariant distribution of the dynamics with exponential ergodicity. Numerical experiments are presented to confirm the theoretical findings and to show the superiority of the algorithm in sampling.
} 
\\
\textbf{AMS subject classifications: } {\rm\small 60H35, 
65C30.}\\

\textbf{Key Words: }{\rm\small} splitting scalar auxiliary variable scheme; explicit scheme; energy-preserving; exponential integrability; strong convergence; long-time weak convergence; sampling
\end{abstract}

\section{Introduction}\label{section:introduction}
In this paper we are concerned with the following kinetic Langevin dynamics:
\begin{equation}\label{eq:langevin_sde}
\left\{
\begin{array}{l}
\mathrm{d} 
        v(t)=
            -\kappa \nabla \Phi(u(t)) \mathrm{d} t-\gamma v(t) {\rm{d}} t+\Gamma {\rm{d}} W_{t}, 
\\ 
\mathrm{d} 
        u(t)=
            v(t) {\rm{d}} t,
            \quad
            t>0,
\end{array}\right. 
    \end{equation}
subject to initial conditions $v(0)=v_0$ and $u(0)=u_0$,
where $v_0,u_0 \colon \Omega \rightarrow \mathbb{R}^m, m \in \mathbb{N}$ are two random variables,
%
$\Phi \colon \mathbb{R}^m \rightarrow \mathbb{R}$ is called the potential function,
and
$\kappa, \gamma > 0$, $\Gamma \in \mathbb{R}^{m \times m}$.

Originally derived from physics, the kinetic Langevin dynamics \eqref{eq:langevin_sde} captures the evolution of the position $u(t)$ and momentum $v(t)$ for a particle under a force field $-\kappa \nabla \Phi$, subject to both friction and stochastic forcing.
Specifically, $\Phi(\cdot)$ is called the potential function, $\kappa$ the inverse mass, $\gamma$ the friction coefficient and $\Gamma$ denotes the intensity of the stochastic forcing driven by a standard $m$-dimensional Brownian motion $\{W_{t}\}_{t\geq 0}$. 
Under some mild assumptions (see \cite[Proposition 6.1]{pavliotis2014stochastic}) and by setting $\Gamma=\sqrt{2\kappa\gamma} I_{m\times m}$, the model \eqref{eq:langevin_sde} admits a unique invariant measure given by
$$
\mu_{\infty}(x,y)=
\mu^{(1)}_{\infty}(y)
\mu^{(2)}_{\infty}(x)
    \propto
    \exp\big(- \tfrac{|y|^2}{2\kappa}\big) \exp\big(-\Phi(x) \big),\ x,y\in \mathbb{R}^{m}.
$$
By time-rescaling $\tilde{t}=\gamma t$ and letting $\gamma \rightarrow +\infty$,
the system \eqref{eq:langevin_sde} reduces to the classical overdamped Langevin dynamics (see \cite[Section 6.5.1]{pavliotis2014stochastic}):
\begin{equation}\label{eq:overdamped_langevin}
\mathrm{d} u(t)= -\kappa \nabla \Phi(u(t)) \mathrm{d} t
+\sqrt{2\kappa}{\rm d} W_t,
\
t >0,
\end{equation}
which admits an invariant distribution $\tilde{\mu}_{\infty}(x) \propto \exp(-\Phi(x)), x \in \mathbb{R}^{m}$, the same as the marginal distribution for $u$ of the former system \eqref{eq:langevin_sde}.
Both 
the kinetic Langevin dynamics \eqref{eq:langevin_sde} and its overdamped counterpart \eqref{eq:overdamped_langevin}
have been studied extensively in the literature, widely used for sampling from a target distribution $\pi(x) \propto \exp(-\Phi(x)), x \in \mathbb{R}^{m}$.
%
%
%
%
%
%
%
%
The kinetic (or underdamped) Langevin dynamics \eqref{eq:langevin_sde} is known to converge to equilibrium faster than its overdamped version \eqref{eq:overdamped_langevin}; see, e.g., \cite{cao2023explicit,Eberle2019}.
Moreover,
the kinetic-Langevin-based sampling algorithms (termed as underdamped Langevin Monte Carlo)
also achieve better long-time convergence towards the target distribution, compared to those based on overdamped Langevin dynamics (termed as overdamped Langevin Monte Carlo), particularly in the scenario of high-dimensional sampling; see, e.g., \cite{cheng2018underdamped, dalalyan2020sampling,cheng2018sharp}.
Therefore, it is of great interest to examine the kinetic Langevin dynamics \eqref{eq:langevin_sde} and particularly its numerical discretizations for practical use.
%



In the literature, considerable attention has been paid to numerical approximations of \eqref{eq:langevin_sde} with the globally Lipschitz condition imposed on $\nabla \Phi$; see, e.g.,\cite{abdulle2015long,leroy2024adaptive,cheng2018underdamped,cheng2018sharp,dalalyan2020sampling,schuh2024convergence,leimkuhler2024contraction,leimkuhler2024numerical,gao2022global}, where some authors focused on the tasks of sampling and optimization.
By contrast, much fewer works are devoted to handling the more challenging setting when the potential $\Phi$ exhibits superquadratic growth, i.e., $\nabla \Phi$ is of superlinear growth. Distinct from the first-order SDEs \eqref{eq:overdamped_langevin}, the global monotonicity condition would be missing for the system \eqref{eq:langevin_sde} evolving in the product space $\mathbb{R}^{m} \times \mathbb{R}^{m}$, even with globally monotone drift (i.e., $-\kappa \nabla \Phi$) of superlinear growth. This unavoidably causes essential difficulties in strong and weak convergence analysis of numerical schemes (see, e.g., \cite{Hutzenthaler2020} for more explanation).

In the year of $2002$, the authors of pioneering works \cite{talay2002stochastic,mattingly2002ergodicity}
 investigated long-time behaviors of approximating invariant measures of \eqref{eq:langevin_sde} via implicit schemes in a non-globally Lipschitz setting. 
%
For the strong approximations over a finite time horizon,
the authors of \cite{hutzenthaler2018exponential}  proposed a class of stopped increment-tamed Euler–Maruyama methods for general SDEs without globally monotone coefficients, covering the kinetic Langevin dynamics \eqref{eq:langevin_sde} as a special case. 
The proposed schemes are explicit and  possess an exponential integrability property, which essentially enables the authors of \cite{Hutzenthaler2020,dai2025perturbation} to obtain finite time strong convergence rates.
In \cite{cui2022density}, the researchers constructed an implicit scheme, called the splitting averaged vector field (AVF) method,  for the kinetic Langevin
dynamics \eqref{eq:langevin_sde}. The splitting AVF method preserves the exponential integrability and its strong convergence rate of order one was successfully revealed.
Using a different splitting strategy, the article \cite{chen2025new} recently proposed a new implicit scheme  having finite time strong convergence and long-time weak convergence rates. 
When the dynamics \eqref{eq:langevin_sde} was used for sampling and optimization, the authors of \cite{johnston2024kinetic} introduced explicit tamed Euler schemes and obtained a convergence property (depending on $\gamma$) between the distribution of the schemes and the exact invariant measure in $2$-Wasserstein distance, in a setting of non-globally Lipschitz but strongly convex $\nabla \Phi$.

In this work, we 
propose a novel explicit time-stepping scheme
for the kinetic Langevin dynamics \eqref{eq:langevin_sde} with both quadratic and possibly superquadratic potentials in a non-convex setting. 
The proposed method, called splitting scalar auxiliary variable (SSAV) scheme, achieves a strong convergence rate of order-one over finite time intervals and  exhibits long-time weak convergence with respect to the exact invariant measure.
Compared to the existing schemes, the newly proposed one enjoys the advantages of being explicit, computationally efficient, and exhibiting good long-time performance under dissipative-type conditions on $\nabla \Phi$. 
In particular, the novel scheme has a significant advantage when used to sample from a high-dimensional target distribution, where implicit methods are too expensive.
%
%
%
%
We would like to mention that, the scalar auxiliary variable (SAV) method originally introduced by \cite{shen2018scalar,shen2018convergence},
is efficient and energy stable in solving deterministic gradient flows. Recently, its stochastic extensions have been investigated for strongly approximating stochastic partial differential equations (see, e.g., \cite{cui2025stochastic,metzger2024convergent,metzger2025strong}). 
Different from these existing works, the current study is concerned with the kinetic Langevin
dynamics and a new auxiliary variable is introduced (see Remark \ref{rem:new-auxilary-variables} for details), which combined with the splitting strategy results in a new energy-stable splitting scalar auxiliary variable (SSAV) method. Furthermore, the weak convergence of this type of method is examined and weak error estimates with order one in a large time regime are obtained (see Theorem \ref{thm:longtime_weak_conver}):
\begin{equation}
    \begin{aligned}
    &\Big|
            \mathbb{E}\big[
            \varphi(Y^{x_0}_{t_N})
                      \big]
            -
            \int_{\mathbb{R}^{2m}} \varphi {\rm d}\mu_{\infty}
    \Big|
&\leq
C_1
e^{-\lambda t_N}
+
C_2(1+(t_N)^{l})h,
    \end{aligned}
\end{equation} 
where $h$ is the stepsize and the error constant $C_2(1+(t_N)^{l}) $ polynomially (not exponentially) depends on the time length $t_N = Nh$. This is crucial to significantly reducing the computational costs in approximations of the integral with respect to invariant measures (see subsection \ref{subsec:computational-costs}).
%
%
%
%
%
%
%



The remainder of this paper is structured as follows. In the forthcoming section, we introduce the proposed method. In Section \ref{section:energy_exponential_integrability}, we show the energy-preserving and exponential integrability properties of the proposed scheme. Equipped with the exponential integrability, we derive the order-one strong convergence in Section \ref{section:order_one_strong}. In Section \ref{section:long_time_weak}, 
the long-time weak convergence of the scheme is obtained. Moreover, we compare the computational costs of our method with those of existing methods in approximating the invariant measure. The final section provides some numerical experiments to demonstrate the superiority of the SSAV method for sampling.

\section{Preliminaries and the proposed scheme}\label{section:preliminaries_proposed_scheme}

In this section, we begin with some definitions and notations.
%
Given $m\in \mathbb{N}$, let $\{W_t\}_{t\geq 0}$ be an $m$-dimensional 
standard Wiener process defined on the completed probability space $(\Omega,\mathcal{F}, \{\mathcal{F}_t\}_{t\geq 0},\mathbb{P})$  with filtration $\{\mathcal{F}_t\}_{t\geq 0}$. 
By $|\cdot|$ and $\langle \cdot , \cdot \rangle$ we denote the Euclidean norm 
and the inner product on $\mathbb{R}^m$, respectively. 
Let $\| A \|:= \sqrt{\operatorname{tr}(A^{*}A)}$ be the trace norm of a matrix $A \in \mathbb{R}^{r \times q}, r,  q \in \mathbb{N}$, where $A^*$ denotes the transpose of $A$. 
For a random variable $\xi: \Omega \rightarrow \mathbb{R}^z, z\in \mathbb{N}$, $\mathbb{E}[\xi]$ denotes its expectation and for any $p>0$, $\|\xi\|_{L^p(\Omega;\mathbb{R}^z)}:= (\mathbb{E}[|\xi|^p])^{1/p}$.
Given $h\in (0,1]$, we construct a uniform mesh $\{ t_n \}_{n \in \mathbb{N}_0 }$
and denote
\begin{equation}
    \lfloor t \rfloor_h:=\sup_{n \in \mathbb{N}_0 }\{t_n:t_n \leq t\},
    \
    \mathbb{N}_0 : = \mathbb{N} \cup \{0\}.
\end{equation}
Throughout the paper, we define 
$$
\mathbf{0}:= 0 \cdot I_{m\times m},\ 
\overset{\rightarrow}{0}:=
0 \cdot I_{m\times 1}
$$
and use $C$ to denote a generic positive constant independent of $h$ that may change at different appearances.

%




The considered equation \eqref{eq:langevin_sde} is usually treated as a stochastic Hamiltonian system (see, e.g.,  \cite{talay2002stochastic}), with the Hamiltonian function $H$, also called the total energy, defined by
\begin{equation} \label{eq:hamiltonian_quantity}
    H(v(t),u(t)):=\tfrac{|v(t)|^2}{2}+\kappa\Phi(u(t))+C_H,
    \ 
    t \geq 0,
\end{equation}
where $C_H$ is a constant such that $H(v(t),u(t)) \geq 0$ for all $t\geq 0$.  
Unlike the situation for the deterministic Hamiltonian system, the total energy in the stochastic system is not conserved, as established in Lemma \ref{lem:energy_law}. To show this, we make an assumption as follows.

\begin{assumption}\label{ass:strong_solution}
    Assume $\Phi \in  C^2(\mathbb{R}^m;\mathbb{R})$ and there 
    exists a unique strong solution for the kinetic Langevin dynamics
    \eqref{eq:langevin_sde}.
\end{assumption}
\begin{lemma}[Evolution of the total energy]
\label{lem:energy_law} 
Let Assumption \ref{ass:strong_solution} hold and
let $ \{ ( v(t), u(t) ) \}_{t \geq 0}$ be the solution of
the kinetic Langevin dynamics \eqref{eq:langevin_sde}. Then
    \begin{equation}\label{eq:energy_law_inequal} 
        \mathbb{E}\big[H(v(t),u(t))\big] \leq \mathbb{E}\big[H(v_0,u_0)\big]+\tfrac{1}{2}\|\Gamma\|^2t,
        \
        \forall \ t \geq 0.
    \end{equation}
\end{lemma}
\textbf{Proof:}  Thanks to Itô's formula, one deduces
    \begin{equation}
    \begin{aligned}
     &H(v(t),u(t))-H(v_0,u_0)
\\
     &\quad =\int_0^t \big(v(s)^*,\kappa \nabla \Phi(u(s))^* \big)
        \Big(\begin{array}{c}
            -\kappa\nabla \Phi(u(s)) -\gamma v(s)
            \\
            v(s)
            \end{array}\Big) {\rm d}s
        +
        \int_0^t \big(v(s)^*,\kappa\nabla \Phi(u(s))^* \big)
        \Big( 
        \begin{array}{c}
        \Gamma 
        \\
        \mathbf{0}
        \end{array} 
            \Big) {\rm d}{W}_s
\\
    &\qquad +
    \tfrac{1}{2}\int_0^t 
    {\rm tr}\bigg(
        \Big(\begin{array}{c}
            \Gamma \Gamma^* \  \mathbf{0}
            \\
            \mathbf{0}\ \ \ \  \mathbf{0}
            \end{array}\Big)
             \Big(\begin{array}{c}
            \hspace{-1.5cm} I_{m\times m} \quad \quad \quad   \mathbf{0}
            \\
             \quad \mathbf{0}\quad \quad \kappa{\rm Hess}_x(\Phi(u(s)))
            \end{array}\Big)
    \bigg)
     {\rm d}s
\\
    &\quad=
    \int_0^t
    \Big(
    -\gamma |v(s)|^2 + \tfrac{1}{2}\|\Gamma\|^2
    \Big)
    {\rm d}s
    +
    \int_0^t  
        v(s)^*\Gamma 
{\rm d}{W}_s
\\
&\quad \leq 
\int_0^t  
        v(s)^*\Gamma 
{\rm d}{W}_s
+\tfrac{1}{2} \|\Gamma\|^2t. 
    \end{aligned}
    \end{equation} 
By combining a stopping time argument with Fatou's lemma,  we take expectations on both sides to obtain the desired result.  \qed

To approximate the equation \eqref{eq:langevin_sde}, we first decompose the original equation into a deterministic Hamiltonian system and an Ornstein--Uhlenbeck process, namely, for $t >0$,
\begin{equation}\label{eq:split_origin_langevin}
    \begin{aligned}
            {\rm d}
            \Big(
            \begin{array}{c}
            v(t)
            \\
            u(t)
        \end{array}\Big)
        =
            \Big(
            \begin{array}{c}
            -\kappa\nabla \Phi(u(s))
            \\
            v(t)
        \end{array}\Big){\rm d}t
        +
        \Big(
            \begin{array}{c}
            -\gamma v(t){\rm d}t + \Gamma {\rm d}W_t
            \\
            \overset{\rightarrow}{0}
        \end{array}\Big).
    \end{aligned}
\end{equation}
In order to construct the SSAV scheme, we further make the following assumption.
\begin{assumption}\label{ass:well_pose_rho_t}
For sufficiently large $C_H > 0$, 
there exists some constant $\alpha>0$ such that
\begin{equation}
\kappa \Phi(x)+ C_H- \alpha |x|^{2}  \geq 1, \ \forall x\in \mathbb{R}^m.
\end{equation}
\end{assumption}
We mention that a similar assumption can be found in \cite[Assumption 2.5]{cui2022density}  and \cite[Hypothesis 1.1 (1.19)]{talay2002stochastic}. 
Under this assumption, a key idea is to introduce a new scalar auxiliary variable
\begin{equation} \label{eq:auxiliary_variable}
    \rho(t):= \sqrt{\kappa\Phi(u(t))+C_H- \alpha |u(t)|^2},
\end{equation}
where by assumption we know that $ \rho(t)\geq 1$ is well-defined for all $t\geq 0$.
Then one may rewrite the total energy \eqref{eq:hamiltonian_quantity} as
\begin{equation}\label{eq:energy(withGamma)_introduction}
H(v(t),u(t))
=
\mathcal{H}(v(t),u(t),\rho(t))
:=
\tfrac{|v(t)|^2}{2}+\alpha |u(t)|^2+\rho(t)^2,
\end{equation}
and the system \eqref{eq:split_origin_langevin} will be turned into \begin{equation}\label{eq:split_SAV_origin_equation_fixed_Upsilon}
    \begin{aligned}
            {\rm d}
            \left(
            \begin{array}{c}
            v(t)
            \\
            u(t)
            \\
            \rho(t)
        \end{array}\right)
        =
            \left(
            \begin{array}{c}
            -2\alpha u(t)
           + 
            \tfrac{-\kappa\nabla \Phi(u(t)) +2\alpha u(t)}{\sqrt{\kappa\Phi(u(t))+C_H-\alpha |u(t)|^{2}}}\rho(t)
            \\
            v(t)
            \\
        \Big\langle \tfrac
        {\kappa \nabla \Phi(u(t))-2\alpha u(t)}
        {2\sqrt{\kappa \Phi(u(t))+C_H- \alpha |u(t)|^2}},
        v(t)
        \Big\rangle
        \end{array}\right){\rm d}t
        +
        \left(
            \begin{array}{c}
            -\gamma v(t){\rm d}t + \Gamma {\rm d}W_t
            \\
            \overset{\rightarrow}{0}
            \\
            0
        \end{array}\right).
    \end{aligned}
\end{equation}

Next we
attempt
to develop a novel numerical scheme that numerically preserves the total energy $\mathcal{H}(v(t),u(t),\rho(t))$ of the system \eqref{eq:split_SAV_origin_equation_fixed_Upsilon}.
Given ${v}_n, {u}_n$ and ${\rho}_n$, the first step involves solving
\begin{equation}\label{eq:SSAV_unsolve_equation}
 \left\{\begin{array}{l}
{v}^{\vartriangle}_{n}={v}_{n} 
    -\alpha({u}_{n}+{u}^{\vartriangle}_{n})h 
    +
    \big(
        \tfrac{-\kappa \nabla \Phi({u}_n)+2\alpha {u}_n}{\sqrt{\kappa \Phi({u}_n)+C_H- \alpha |{u}_n|^{2}}} \tfrac{{\rho}^{\vartriangle}_{n}+{\rho}_n}{2}
    \big)h,
\\
 {u}^{\vartriangle}_{n}={u}_n
        +
        \tfrac{{v}^{\vartriangle}_{n}+{v}_n}{2}h,
\\
{\rho}^{\vartriangle}_{n}={\rho}_n
    +
    \tfrac{1}{2}\big\langle \tfrac{\kappa \nabla \Phi({u}_n)-2\alpha {u}_n}{\sqrt{\kappa \Phi({u}_n)+C_H- \alpha |{u}_n|^{2}}},\tfrac{{v}^{\vartriangle}_{n}+{v}_n}{2}\big\rangle 
    h,
    \end{array}\right.
\end{equation}
and the approximations $({v}_{n+1}^*,{u}_{n+1}^*,{\rho}_{n+1})^*$ for the second step are determined by
\begin{equation}\label{eq:SSAV_OU_equation}
\left\{\begin{array}{l}
{v}_{n+1}=
    e^{-\gamma h} {v}^{\vartriangle}_{n}
    +
    \displaystyle{
    \int_{t_n}^{t_{n+1}}} e^{-\gamma(t_{n+1}-s)}\Gamma {\rm{d}} W_{s},
\\
 {u}_{n+1}={u}^{\vartriangle}_{n},
\\
{\rho}_{n+1}={\rho}^{\vartriangle}_n,
    \end{array}\right.
\end{equation}
with the initial values ${v}_0, {u}_0$ and ${\rho}_0=\sqrt{\kappa \Phi(u_0)+C_H-\alpha |u_0|^2}$.
After a glance at \eqref{eq:SSAV_unsolve_equation}, the proposed SSAV scheme \eqref{eq:SSAV_unsolve_equation} appears to be implicit.
In the next proposition, we demonstrate that the  SSAV scheme \eqref{eq:SSAV_unsolve_equation}-\eqref{eq:SSAV_OU_equation} is, in fact, an explicit one. 
%
\begin{proposition}\label{propo:explicit_solvable_SAV}
Under Assumption \ref{ass:well_pose_rho_t},
the equation \eqref{eq:SSAV_unsolve_equation} can be explicitly solvable as follows: 
\begin{equation}\label{eq:solution_split_SAV}
    \left\{\begin{array}{l}
    \rho^{\vartriangle}_{n}=\rho_n +\tfrac{\langle Q_n, 2v_n-2\alpha u_nh -2Q_n\rho_n h\rangle }{2+\alpha h^2+|Q_n|^2h^2}h,
\\ 
    v_n^{\vartriangle} = 
    \tfrac{2v_n-4Q_n I^h_n h
            -
            4\alpha u_nh{}
            -
            \alpha v_n h^2
            }
    {2+\alpha h^2},
\\
    u^{\vartriangle}_{n}=u_n+
    \tfrac{2v_n-2Q_n I_n^h h
            -
            2\alpha u_nh
            }
    {2+\alpha h^2}h
    ,
    \end{array}\right.
    \end{equation}
where, for $n \in \mathbb{N}_0$ we denote 
\begin{align}
\label{eq:Q_n_definition}
Q_n
& :=
\tfrac{\kappa \nabla \Phi(u_n)-2\alpha u_n}{2\sqrt{\kappa \Phi(u_n)+C_H- \alpha |u_n|^{2}}},
\\
\label{eq:I_n^h_definition} 
    I^h_n
    &
    :=
    \tfrac{\rho_n(2+\alpha h^2)+
                        \langle 
                            Q_n, v_n-\alpha u_nh
                        \rangle h}
            {2+\alpha h^2+|Q_n|^2 h^2}.
\end{align}
\end{proposition}

\textbf{Proof:}
First of all, we rewrite the equation \eqref{eq:SSAV_unsolve_equation}
as follows: 
\begin{equation}
\label{eq:determin_part_SAV_scheme_simple}
    \left\{\begin{array}{l}
v^{\vartriangle}_{n}=v_n 
    -\alpha(u_n+u^{\vartriangle}_{n})h 
    -
        Q_n(\rho^{\vartriangle}_{n}+\rho_n) h,
\\
 u^{\vartriangle}_{n}=u_n
        +
        \tfrac{v^{\vartriangle}_{n}+v_n}{2}h,
\\
\rho^{\vartriangle}_{n}=\rho_n
    +\tfrac{1}{2}
   \langle Q_n, {v^{\vartriangle}_{n}+v_n}\rangle 
    h.
    \end{array}\right.
\end{equation}
In view of the first and the second equalities in \eqref{eq:determin_part_SAV_scheme_simple}, we get
\begin{equation}
    v^{\vartriangle}_{n}=v_n 
    -Q_n(\rho^{\vartriangle}_{n}+\rho_n) h
    -\alpha h (2u_n+\tfrac{v^{\vartriangle}_{n}+v_n}{2}h),
\end{equation} or, equivalently,
\begin{equation}\label{eq:proof_solvable_split_v_star}
\begin{aligned}
    v^{\vartriangle}_{n}
=
    \tfrac{2v_n 
    -2Q_n(\rho^{\vartriangle}_{n}+\rho_n) h-4\alpha u_nh-\alpha v_n h^2}
    {2+\alpha h^2}
=v_n-
\tfrac{
    2Q_n(\rho^{\vartriangle}_{n}+\rho_n) h+4\alpha u_nh+2\alpha v_n h^2}
    {2+\alpha h^2}.
\end{aligned}
\end{equation}
This, together with the last equality in \eqref{eq:determin_part_SAV_scheme_simple}, implies that
%
\begin{equation}
\rho^{\vartriangle}_{n}=\rho_n+
h
\Big\langle
Q_n,\tfrac{2v_n-Q_n(\rho^{\vartriangle}_{n}+\rho_n)h-2\alpha u_n h}{2+\alpha h^2}
\Big\rangle.
\end{equation}
Since $\rho_n, \rho^{\vartriangle}_{n}$ are both scalar variables,
the solution of 
$\rho^{\vartriangle}_{n}$
can be easily solved:
\begin{equation}\label{eq:express_rho_n_add_1}
\begin{aligned}
\rho^{\vartriangle}_{n} 
=
\tfrac{(2+\alpha h^2)\rho_n
                       +
                       \langle 
                            Q_n, 2v_n-2\alpha u_nh-Q_n \rho_n h
                        \rangle h
                             }
                    {2+\alpha h^2+|Q_n|^2h^2}
=
\rho_n +\tfrac{\langle Q_n, 2v_n-2\alpha u_nh -2Q_n\rho_n h\rangle }{2+\alpha h^2+|Q_n|^2h^2}h.
\end{aligned}
\end{equation}
Recalling \eqref{eq:I_n^h_definition} and the equation \eqref{eq:proof_solvable_split_v_star} gives the solution for $v^{\vartriangle}_{n}$ that
\begin{equation}
v_n^{\vartriangle} = 
    \tfrac{2v_n-4Q_n I^h_n h
            -
            4\alpha u_nh{}
            -
            \alpha v_n h^2
            }
    {2+\alpha h^2}.
\end{equation}
Plugging this into the second equality of \eqref{eq:determin_part_SAV_scheme_simple}
gives the explicit expression for $u^{\vartriangle}_{n}$.
\qed

By virtue of Proposition \ref{propo:explicit_solvable_SAV}, \eqref{eq:determin_part_SAV_scheme_simple} and 
\eqref{eq:proof_solvable_split_v_star}, one can recast the SSAV scheme as follows: 
\begin{equation}\label{eq:practical_iteration_SSAV}
    \left\{\begin{array}{l}
    \rho_{n+1}=\rho_n +\tfrac{\langle Q_n, 2v_n-2\alpha u_nh -2Q_n\rho_n h\rangle }{2+\alpha h^2+|Q_n|^2h^2}h,
\\ 
    v_n^{\vartriangle} =
    v_n-
    \tfrac{
    2Q_n(\rho_{n+1}+\rho_n) h+4\alpha u_nh+2\alpha v_n h^2}
    {2+\alpha h^2},
\\ 
    v_{n+1} = 
    e^{-\gamma h}
    v_n^{\vartriangle}
    +
    \sqrt{
    \tfrac{1}{2\gamma}
    (
    1-e^{-2\gamma h}
    )
    }
    \Gamma
    \mathbbm{N}_m
    ,
\\
    u_{n+1}=u_n+
   \tfrac{v_n^{\vartriangle}+v_n}{2}h,
    \end{array}\right.
\end{equation}
where 
$\mathbbm{N}_m \sim \mathbf{N}(\overset{\rightarrow}{0}, I_{m\times m})$ and $Q_n$ is defined by \eqref{eq:Q_n_definition} for $n \in \mathbb{N}_0$. This way one can only calculate $Q_n$ for different potential functions, facilitating the implementation of the proposed scheme.


Before closing this section, we add some necessary remarks on the scalar auxiliary variable $\rho$.
\begin{remark}
\label{rem:new-auxilary-variables}
    Different from usual auxiliary variables like $\tilde{\rho}(t):=\sqrt{\kappa \Phi(u(t))+ C_H}, t\geq 0$ introduced in existing works \cite{shen2018scalar,cui2025stochastic,metzger2024convergent}, by subtracting $\alpha |u(t)|^2$ we develop a new auxiliary variable $\rho (t) := \sqrt{\kappa\Phi(u(t))+C_H- \alpha |u(t)|^2}$, given by \eqref{eq:auxiliary_variable}. 
Thanks to the introduce of $\alpha |u(t)|^2$, one can derive the exponential integrability property of the numerical approximation $u_n$ (cf. Theorem \ref{thm:exponential_integ_numer_solution}), which is essentially required in recovering the strong convergence rate. Otherwise, only the exponential integrability property of $v_n$ is available.
Also, it should be noted that, the auxiliary variable $\rho$ and our approach are still new, even for the deterministic Hamiltonian setting, i.e., $\gamma=\Gamma=0$.
\end{remark}


\section{Energy-preserving and exponential integrability}\label{section:energy_exponential_integrability}

In this section, we attempt to show the energy-preserving and exponential integrability properties of the proposed scheme.
%
%

\subsection{Energy-preserving}



The next theorem concerns
the energy-preserving property of the iteration \eqref{eq:SSAV_unsolve_equation}.

\begin{theorem}\label{thm: numer_energy_preserv}
Under Assumption \ref{ass:well_pose_rho_t}, the iteration \eqref{eq:SSAV_unsolve_equation} preserves the modified energy, that is,
\begin{equation}\label{eq:energ_validat_0}
\tfrac{|v^{\vartriangle}_{n}|^2}{2}+\alpha |u^{\vartriangle}_{n}|^2+\rho^{\vartriangle2}_{n}=
    \tfrac{|v_{n}|^2}{2}+\alpha |u_{n}|^2+\rho_{n}^2,
\quad
\forall \ n \in \mathbb{N}_0.
\end{equation}
\end{theorem}
\textbf{Proof:} 
Since the identity \eqref{eq:energ_validat_0} can be rearranged as
\begin{equation}\label{eq:energ_validat}
\tfrac{1}{2} \langle v^{\vartriangle}_n+v_n, v^{\vartriangle}_n-v_n \rangle
+
\alpha \langle u^{\vartriangle}_{n}+u_n, u^{\vartriangle}_{n}-u_n\rangle
=
\rho_{n}^2-\rho^{\vartriangle2}_{n},
\end{equation}
it remains to validate \eqref{eq:energ_validat}.
In light of \eqref{eq:proof_solvable_split_v_star}, for all $n \in \mathbb{N}_0$ we compute
\begin{equation}
\left\{\begin{array}{l}
    v^{\vartriangle}_n+v_n= \tfrac{4v_n-2(\rho^{\vartriangle}_{n}+\rho_{n})Q_nh-4\alpha u_nh}{2+\alpha h^2},
\\ 
    v_n^{\vartriangle} -v_n= 
    \tfrac{-2Q_n(\rho^{\vartriangle}_{n}+\rho_{n}) h
            -
            4\alpha u_nh
            -
            2\alpha v_n h^2
            }
    {2+\alpha h^2},
\\ 
    u^{\vartriangle}_{n}+u_{n}
    =
    \tfrac{2v_nh-Q_n(\rho^{\vartriangle}_{n}+\rho_{n})h^2+4u_n
            }
    {2+\alpha h^2},
\\
    u^{\vartriangle}_{n}-u_n=
    \tfrac{2v_nh-
        Q_n(\rho^{\vartriangle}_{n}+\rho_{n})h^2
            -
            2\alpha u_nh^2
            }
    {2+\alpha h^2}.
    \end{array}\right.
\end{equation}
Therefore, one can handle the left side of \eqref{eq:energ_validat} in the following way:
\begin{equation}
\begin{aligned}
    &\tfrac{1}{2} \langle v^{\vartriangle}_n+v_n, v^{\vartriangle}_n-v_n \rangle
    +
    \alpha \langle u^{\vartriangle}_{n}+u_n, u^{\vartriangle}_{n}-u_n\rangle
\\
    &\quad =
    \tfrac{1}{(2+\alpha h^2)^2}
    \Big(
    \left\langle 2v_nh-(\rho^{\vartriangle}_{n}+\rho_{n})Q_nh^2-2\alpha u_nh^2
    , 
            -2Q_n(\rho^{\vartriangle}_{n}+\rho_{n}) 
            -
            4\alpha u_n
            -
            2\alpha v_n h
    \right\rangle
\\
    &\qquad
        \qquad
        \qquad
    +
    \left\langle
    2v_nh-(\rho^{\vartriangle}_{n}+\rho_{n})Q_nh^2-2\alpha u_nh^2
    ,
    2v_n\alpha h-\alpha Q_n(\rho^{\vartriangle}_{n}+\rho_{n})h^2+4\alpha u_n
    \right\rangle
    \Big)
\\
&\quad=
\tfrac{1}{(2+\alpha h^2)^2}
    \left\langle 2v_nh-(\rho^{\vartriangle}_{n}+\rho_{n})Q_nh^2-2\alpha u_nh^2
    , 
            -2Q_n(\rho^{\vartriangle}_{n}+\rho_{n}) 
            -\alpha Q_n(\rho^{\vartriangle}_{n}+\rho_{n})h^2
    \right\rangle 
\\
&\quad =
    - \tfrac{\rho^{\vartriangle}_{n}+\rho_{n}}{2+\alpha h^2}
    \left\langle 2v_nh-(\rho^{\vartriangle}_{n}+\rho_{n})Q_nh^2-2\alpha u_nh^2
    , 
     Q_n
    \right\rangle,
\end{aligned}
\end{equation}
where,
by \eqref{eq:express_rho_n_add_1} we arrive at
\begin{equation}
\begin{aligned}
    &
    - \tfrac{1}{2+\alpha h^2}
    \left\langle 
    2v_nh-(\rho^{\vartriangle}_{n}+\rho_{n})Q_nh^2-2\alpha u_nh^2
    , 
     Q_n
    \right\rangle
\\
    &\quad =
    - \tfrac{1}{2+\alpha h^2}
    \left\langle 2v_n h - 2\alpha u_n h^2
    , 
     Q_n
    \right\rangle
    +
    \left(
    \tfrac{|Q_n|^2h^2}{2+\alpha h^2} 
    + 1
    \right)
    (\rho^{\vartriangle}_{n}+\rho_{n})
    -
    (\rho^{\vartriangle}_{n}+\rho_{n})
\\
    &\quad = \rho_{n}-\rho^{\vartriangle}_{n},
\end{aligned}
\end{equation}
and the assertion \eqref{eq:energ_validat} follows.
\qed

As a by-product of Theorem \ref{thm: numer_energy_preserv}, one has the following evolution of the numerical total energy.

\begin{corollary}[Evolution of the numerical total  energy]\label{corol:num_energy_law}
Let Assumption \ref{ass:well_pose_rho_t} hold
and let
the SSAV scheme
$\{({v}_n^*,{u}_n^*,{\rho}_n)^*\}_{ n \in \mathbb{N}_0 }$
be produced by \eqref{eq:SSAV_unsolve_equation}-\eqref{eq:SSAV_OU_equation}.
Then
\begin{equation}
    \mathbb{E}\big[
    \mathcal{H}({v}_{n},{u}_{n},{\rho}_{n})
    \big]
    \leq 
    \mathbb{E}\big[
    \mathcal{H}(v_0,u_0,\rho_0)
    \big]
    +
    \tfrac{1}{2}\|\Gamma\|^2t_{n},
    \quad
    \forall \ n \in \mathbb{N}_0.
\end{equation}
\end{corollary}
\textbf{Proof:} For fixed $n \in \mathbb{N}_0$ and $s\in [t_n,t_{n+1}]$, we define
\begin{equation}\label{eq:def_v_s_n}
{v}_{s,n}:=
v_n^{\vartriangle}
+
\int_{t_n}^{s}
-\gamma {v}_{r,n}
{\rm d}r
+
\int_{t_n}^{s}
\Gamma  
{\rm d}W_r.
\end{equation}
Given $v_n^{\vartriangle}$, one clearly observes
${v}_{n+1}={v}_{t_{n+1},n}$
from
the first equality in \eqref{eq:SSAV_OU_equation}.
%
Applying Theorem \ref{thm: numer_energy_preserv} and using the It\^{o} formula, one deduces
\begin{equation}\label{eq:evolution_energy_proof_SSAVA1}
    \begin{aligned}
    &\mathbb{E} 
        \big[
        \tfrac{|{v}_{n+1}|^2}{2}
        +
        \alpha |{u}_{n+1}|^2
        +
        {\rho}_{n+1}^2
        \big]
\\
    &\quad =\mathbb{E} 
        \big[
        \tfrac{|v^{\vartriangle}_{n}|^2}{2}
        +
        \alpha |u^{\vartriangle}_{n}|^2
        +
        \rho^{\vartriangle2}_{n}
        \big]
        +
        \mathbb{E} 
        \big[
        \tfrac{|{v}_{n+1}|^2}{2}
        -
        \tfrac{|v^{\vartriangle}_{n}|^2}{2}
        \big]
\\
    &\quad =\mathbb{E} 
        \big[
        \tfrac{|{v}_{n}|^2}{2}
        +
        \alpha |{u}_{n}|^2
        +
        {\rho}_{n}^2
        \big]
        +
        \mathbb{E} 
        \Big[
        \int_{t_n}^{t_{n+1}}
            \big(
            -\gamma |{v}_{s,n}|^2
            +\tfrac{1}{2}\|\Gamma\|^2
            \big)
        {\rm d}s
        \Big]
        +
        \mathbb{E} 
        \Big[
        \int_{t_n}^{t_{n+1}}
            {v}_{s,n}^* \Gamma
        {\rm d}W_s
        \Big]
\\
&\quad \leq 
\mathbb{E} 
        \big[
        \tfrac{|{v}_{n}|^2}{2}
        +
        \alpha |{u}_{n}|^2
        +
        {\rho}_{n}^2
        \big]
+
\tfrac{1}{2} \|\Gamma\|^2 h
\\
&\quad \leq 
\mathbb{E} 
        \big[
        \tfrac{|v_{0}|^2}{2}
        +
        \alpha |u_{0}|^2
        +
        \rho_{0}^2
        \big]
+
\tfrac{1}{2} \|\Gamma\|^2 t_{n+1},
    \end{aligned}
\end{equation}
which completes the proof.\qed 

\begin{remark}
    We emphasise that, in general, the numerical total energy $\mathcal{H}(v_n,u_n,\rho_n)$ is not equal to the exact total energy $\mathcal{H}(v(t_n),u(t_n),\rho(t_n))$ (their discrepancies are bounded in Proposition \ref{propo:energy_error}).
\end{remark}

\subsection{Exponential integrability properties}

In the following, for any $T\in(0,+\infty)$ we construct a uniform mesh on the finite time horizon $[0, T]$, with the uniform stepsize $h := T/ N_T$ for any $N_T \in \mathbb{N}$.
This subsection aims to establish exponential integrability properties of both the exact and numerical solutions.
To this end, we present an important lemma following from \cite[Corollary 2.4]{cox2024local}.
\begin{lemma}\label{lem:exponen_integ_lemma}
    Assume that $0\leq t_0<T , \lambda \in \mathbb{R},U_0 \in C^2( \mathbb{R}^d, \mathbb{R})$ and $U_1 \in C( \mathbb{R}^d, \mathbb{R})$ for $d \in \mathbb{N}$. 
    Let $\{Z(t)\}_{t\in[t_0,T]}$ be an adapted stochastic process with continuous sample paths fulfilling 
    $
    \int_{t_0}^T |\mathbb{F}(s,Z(s))|+\|\mathbb{G}(s,Z(s))\|^2  {\rm d}s < +\infty \ \mathbb{P}$-a.s.,
    and for any  $t\in[t_0,T]$,
    $$
    Z(t)=Z(t_0)+  
    \int_{t_0}^{t} \mathbb{F}(s,Z(s)) {\rm d}s
    +
    \int_{t_0}^{t} \mathbb{G}(s,Z(s)) {\rm d}W_s,
    \ \mathbb{P}\text{-a.s.}.
    $$
    In addition, suppose that for any $x\in \mathbb{R}^d, t\in [t_0,T]$,
    \begin{equation}\label{eq:exact_exponen_integ_condition}
     U_0'(x) \mathbb{F}(t,x) 
     +
     \tfrac{1}{2} \operatorname{tr}
     \left(\mathbb{G}(t,x)\mathbb{G}(t,x)^{*}\operatorname{Hess}_x\big(U_0(x) \big)\right)
     +
     \tfrac{1}{2e^{\lambda t}}|\mathbb{G}(t,x)^* \nabla U_0(x)|^2
    +U_1(x)
    \leq 
    \lambda U_0(x).
     \end{equation}
    Then $\{Z(t)\}_{t\in[t_0,T]}$ possesses the property
    \begin{equation}
    \sup_{t\in[t_0,T]}
    \mathbb{E} 
        \Big[\exp
        \big(
            e^{-\lambda t}U_0(Z(t))+ \int_{t_0}^t e^{-\lambda s}U_1(Z(s)) {\rm d}s
        \big) 
        \Big]
    \leq 
    \mathbb{E} 
        \Big[\exp\big(e^{-\lambda t_0}{U_0(Z({t_0}}))\big)
        \Big].
    \end{equation}
\end{lemma}
\qed

Equipped with Lemma \ref{lem:exponen_integ_lemma}, one can obtain the exponential integrability property of the solution of the kinetic Langevin dynamics \eqref{eq:langevin_sde}, as stated in the forthcoming lemma.

\begin{lemma}\label{lem:exact_solution_exponen_property}
Let  Assumptions \ref{ass:strong_solution}, \ref{ass:well_pose_rho_t} hold
and assume 
$\mathbb{E}\big[
        \exp(\delta(\tfrac{|v_0|^2}{2} 
+\kappa \Phi(u_0)+C_H))
        \big]<+\infty$
for some $\delta>0$,
Then for any $\lambda \geq 
\max\{ \delta \|\Gamma\|^2-2\gamma, 0\}$, the solution of the kinetic Langevin dynamics \eqref{eq:langevin_sde} obeys
\begin{equation}\label{eq:expon_integ_prope_exact_solution}
\begin{aligned}
&\sup_{t\in[0,T]}
\mathbb{E}
\Big[
        \exp\big(
            e^{-\lambda t}
            \delta(\tfrac{|v(t)|^2}{2} + \alpha |u(t)|^2+\rho(t)^2)
            \big)
\Big]
\\
&\quad \leq 
\exp
        \Big(
        \tfrac{\delta}{2}\|\Gamma\|^2
        \big(
        \mathbbm{1}_{\{\lambda=0\}}T
        +
        \mathbbm{1}_{\{\lambda > 0 \}}
        \tfrac{1}{\lambda}
        \big)
        \Big)
\mathbb{E}\Big[
        \exp\big({\delta(\tfrac{|v_0|^2}{2} +\kappa \Phi(u_0)+C_H)}\big)
        \Big].
\end{aligned}
\end{equation}
\end{lemma}
\textbf{Proof:} For the case $\lambda>0$, we refer to \cite[Proposition 3.3]{cui2022density},
whereas the proof for the case $\lambda = 0$ follows the same lines as the case $\lambda>0$ and is therefore omitted.
\qed

Next, recalling $N_T\in \mathbb{N}$ and $h=T/N_T$,
we establish the exponential integrability property for the SSAV scheme \eqref{eq:SSAV_unsolve_equation}-\eqref{eq:SSAV_OU_equation}.

\begin{theorem}\label{thm:exponential_integ_numer_solution}
Let Assumption \ref{ass:well_pose_rho_t} hold and
assume
$\mathbb{E}\big[
        \exp(\delta(\tfrac{|v_0|^2}{2} +\kappa \Phi(u_0)+C_H))
        \big]<+\infty$
for some $\delta>0$.
Then for any $\lambda \geq 
\max\{ \delta \|\Gamma\|^2-2\gamma, 0\}$, the SSAV scheme $\{({v}_n^*,{u}_n^*,{\rho}_n)^*\}_{ n \in \mathbb{N}_0}$ defined by \eqref{eq:SSAV_unsolve_equation}-\eqref{eq:SSAV_OU_equation} obeys
%
\begin{equation}\label{eq:expon_integ_prope_numerical_solution}
\begin{aligned}    
&\sup_{h\in(0,1]}\sup_{n\in\{0,1,...,N_T\}}
\mathbb{E}
\Big[
    \exp\big(
        e^{-\lambda t_n}
        \delta(\tfrac{|{v}_n|^2}{2} + \alpha |{u}_n|^2+{\rho}_n^2)
        \big)
\Big]
\\
&\quad \leq 
    \exp
        \Big(
        \tfrac{\delta}{2}\|\Gamma\|^2
        \big(
        \mathbbm{1}_{\{\lambda=0\}}T
        +
        \mathbbm{1}_{\{\lambda > 0 \}}
        \tfrac{1}{\lambda}
        \big)
        \Big)
\mathbb{E}\Big[
        \exp\big({\delta(\tfrac{|v_0|^2}{2} +\kappa \Phi(u_0)+C_H)}\big)
        \Big].
\end{aligned}
\end{equation}
\end{theorem}
\textbf{Proof:} 
For $t\in[t_n,t_{n+1}], n \in \mathbb{N}_0$, we consider the following systems
\begin{equation}\label{eq:numer_linear_system}
\left\{\begin{array}{l}
\begin{aligned}
    &{\rm d} {v}_{t,n} = -\gamma {v}_{t,n} {\rm d}t +\Gamma {\rm d}W_t,
\\
    &{\rm d} {u}_t =  \overset{\rightarrow}{0},
\\
    &{\rm d} {\rho}_t= 0,
\end{aligned}
    \end{array}\right.
\end{equation}
initiating at $v^{\vartriangle}_n,u^{\vartriangle}_{n}$ and $\rho^{\vartriangle}_{n}$.
For $x= (x_1^*,x_2^*,x_3)^*, \ x_1, x_2\in \mathbb{R}^m,x_3\in \mathbb{R}$, we set
$$
U_0(x)=\delta (\tfrac{|x_1|^2}{2}+\alpha |x_2|^2+x_3^2),
\
{U}_1(x)=-\tfrac{1}{2}\delta \|\Gamma\|^2,
$$
and for $t \geq 0$,
$$
{\mathbb{F}}(t,x)= (-\gamma x_1^*, \overset{\rightarrow}{0}\vspace{0.01pt}^*, 0)^*,
\
{\mathbb{G}}(t,x)=(\Gamma^*, \mathbf{0}, \overset{\rightarrow}{0})^*.
$$
Now we check the condition \eqref{eq:exact_exponen_integ_condition}, so that the exponential integrability property of the system \eqref{eq:numer_linear_system} can be obtained due to Lemma \ref{lem:exponen_integ_lemma}. For
$\lambda \geq 
\max\{ \delta \|\Gamma\|^2-2\gamma, 0\}$, it is easy to deduce
\begin{equation}
\begin{aligned}
&
U_0'(x){\mathbb{F}}(x)
    +
    \tfrac{1} {2}\operatorname{tr}\big({\mathbb{G}}(x){\mathbb{G}}(x)^{*}\operatorname{Hess}_x\big(U_0(x)\big)\big)
    +
    \tfrac{1}{2e^{\lambda t}}|{\mathbb{G}}(x)^* \nabla U_0(x)|^2
    +{U}_1(x)
\\
&\quad =
-\delta \gamma
|x_1|^2
+
\tfrac{1}{2}\delta \|\Gamma\|^2
+
\tfrac{1}{2}{e^{-\lambda t}} \delta^2
|\Gamma^*x_1|^2
-
\tfrac{1}{2}\delta \|\Gamma\|^2
\\
&\quad \leq
(
\delta
\|\Gamma\|^2
- 2\gamma
)
\delta
\tfrac{|x_1|^2}{2}
\\
&\quad \leq \lambda U_0(x).
\end{aligned}
\end{equation}
Using Lemma \ref{lem:exponen_integ_lemma} and 
Theorem \ref{thm: numer_energy_preserv} hence gives
\begin{equation}
\begin{aligned}
 &\mathbb{E}\big[
        \exp\big(e^{-\lambda t_{n+1}}\delta(\tfrac{|{v}_{n+1}|^2}{2} + \alpha |{u}_{n+1}|^2
        +
        {\rho}_{n+1}^2)\big)
    \big]
\\
&\quad \leq 
        \mathbb{E}\big[
        \exp\big(e^{-\lambda t_n}
        \delta(\tfrac{|v^{\vartriangle}_n|^2}{2} + \alpha |u^{\vartriangle}_{n}|^2
        +
        \rho^{\vartriangle2}_{n}
        )\big)
        \big]
        \exp
        \Big(
        \tfrac{\delta}{2}\|\Gamma\|^2
        \int_{t_n}^{{t_{n+1}}}
        e^{-\lambda s} {\rm d}s
        \Big)
\\
&\quad =
\mathbb{E}\big[
        \exp\big(e^{-\lambda t_n}
        \delta(\tfrac{|v^{\vartriangle}_n|^2}{2} + \alpha |u^{\vartriangle}_{n}|^2
        +
        \rho^{\vartriangle2}_{n}
        )\big)
        \big]
        \exp
        \Big(
        \tfrac{\delta}{2}\|\Gamma\|^2
        \big(
        \mathbbm{1}_{\{\lambda=0\}}h
        +
        \mathbbm{1}_{\{\lambda > 0 \}}
        \tfrac{1}{\lambda}
        (e^{-\lambda t_n}-e^{-\lambda t_{n+1}})
        \big)
        \Big)
\\
&\quad =
\mathbb{E}\big[
        \exp\big(e^{-\lambda t_n}
        \delta(\tfrac{|v_n|^2}{2} + \alpha |u_{n}|^2
        +
        \rho_{n}^2
        )\big)
        \big]
    \exp
        \Big(
        \tfrac{\delta}{2}\|\Gamma\|^2
        \big(
        \mathbbm{1}_{\{\lambda=0\}}h
        +
        \mathbbm{1}_{\{\lambda> 0 \}}
        \tfrac{1}{\lambda}
        (e^{-\lambda t_n}-e^{-\lambda t_{n+1}})
        \big)
        \Big)
\\
&\quad \leq 
\mathbb{E}\big[
        \exp\big(
        \delta(\tfrac{|v_0|^2}{2} + \alpha |u_{0}|^2
        +
        \rho_{0}^2
        )\big)
        \big]
    \exp
        \Big(
        \tfrac{\delta}{2}\|\Gamma\|^2
        \big(
        \mathbbm{1}_{\{\lambda=0\}}T
        +
        \mathbbm{1}_{\{\lambda > 0 \}}
        \tfrac{1}{\lambda}
        \big)
        \Big).
\end{aligned}
\end{equation}
The proof is thus completed. \qed

By applying Taylor's expansion for the exponential function, we obtain a corollary that follows directly from Lemma \ref{lem:exact_solution_exponen_property} and Theorem \ref{thm:exponential_integ_numer_solution}.

%
%
%
\begin{corollary}\label{coro:p_order_and_expone_moment}
Let Assumptions \ref{ass:strong_solution}, \ref{ass:well_pose_rho_t} be satisfied and
assume 
$\mathbb{E}\big[
        \exp(\delta(\tfrac{|v_0|^2}{2}  +\kappa \Phi(u_0)+C_H))
        \big]<+\infty$
for some $\delta>0$.
Let the processes $\{(v(t)^*, u(t)^*)^*\}_{t \in [0,T]}$ and $\{(v^*_n, u^*_n, \rho_n)^*\}_{n \in \{0,1, ..., N_T\}} $ be produced by \eqref{eq:langevin_sde} and \eqref{eq:SSAV_unsolve_equation}-\eqref{eq:SSAV_OU_equation}, respectively.
Then there exist some positive constants $C_{p,T}$ and $C_{T}$ such that
\begin{equation}
    \begin{aligned}
      \sup_{t \in [0,T]}
      \mathbb{E}\big[|v(t)|^p+|u(t)|^p\big]
       {\ \textstyle \bigvee}
                \sup_{h\in(0,1]}\sup_{n\in\{0,1,...,N_T\}}
         \mathbb{E}\big[|v_n|^p+|u_n|^p+|\rho_n|^p\big]
         \leq C_{p,T},
      \end{aligned}
    \end{equation}
and for any $\eta \leq \alpha \delta e^{-\lambda T}$,
\begin{equation}
    \begin{aligned}
      \sup_{t \in [0,T]}
      \mathbb{E}\big[\exp(\eta|u(t)|^2)\big]
       {\ \textstyle \bigvee}
                \sup_{h\in(0,1]}\sup_{n\in\{0,1,...,N_T\}}
         \mathbb{E}\big[\exp(\eta|{u}_n|^2)\big]
         \leq C_{T}.
      \end{aligned}
    \end{equation}
\end{corollary}
We would like to point out that the constant $C_{p,T}$ in Corollary \ref{coro:p_order_and_expone_moment} is growing exponentially with respect to $T$. Nonetheless, this is sufficient for analyzing strong convergence of the numerical method over a finite time horizon.
In comparison, 
Proposition \ref{prop:bounded_moment_exact_numer_solu} (see Section 5 below) provides 
moment bounds with only polynomial dependence on $T$,
offering a refined estimate
that is essential for reducing the computational costs in the large time regime, as will be discussed later.

\section{Order-one strong convergence in finite time horizon}\label{section:order_one_strong}
This section is devoted to revealing strong convergence rate of the proposed SSAV scheme \eqref{eq:SSAV_unsolve_equation}-\eqref{eq:SSAV_OU_equation} in finite time horizon $t\in[0,T]$. Here and throughout this paper, we denote
\begin{equation}\label{eq:define_X_t_and_Y_n}
    X(t):=(v(t)^*,u(t)^*)^*,\
    t \geq 0,\
    \text{and}\
    {Y}_{t_n}:= ({v}_n^*,{u}_n^*)^*,\
    n \in \mathbb{N}_0.
\end{equation}
Also, we rewrite the kinetic  Langevin dynamics \eqref{eq:langevin_sde} 
in the following compact form
\begin{equation}\label{eq:langevin_sde_newform}
    d X(t)
    =
    F(X(t)) dt
    +
    G dW_t,
    \
    t>0,
\end{equation}
where $F$ and $G$ are defined by
\begin{equation}
F(x):=\Big( 
            \begin{array}{c}
            -\kappa \nabla \Phi(x_2) -\gamma x_1
            \\
            x_1
            \end{array} 
            \Big)
\text{ for }
x=
\Big( 
            \begin{array}{c}
            x_1
            \\
            x_2
            \end{array} 
            \Big)\in \mathbb{R}^{2m}
,
\
G:=\Big( 
\begin{array}{c}
    \Gamma 
     \\
    \mathbf{0}
    \end{array} 
            \Big).
\end{equation}
In view of \eqref{eq:solution_split_SAV} and \eqref{eq:SSAV_OU_equation}, 
we introduce a continuous version of
the SSAV method, given by
\begin{equation}\label{eq:continu_ito_process_SSAV}
    \left\{\begin{array}{l}
    \displaystyle{
\hat{v}_{t}=
    {v}_{n}
    +
    \int_{t_n}^{t} 
    \big(
    \tfrac{-4Q_n I_n^h 
            -
            4\alpha {u}_n
            -
            2\alpha {v}_n h
            }
    {2+\alpha h^2}
    -\gamma \hat{v}_s 
    \big)
    {\rm{d}} s
    +
    \int_{t_n}^{t}  \Gamma {\rm{d}} W_{s}},
\\
 \hat{u}_{t}={u}_n
        +
        \displaystyle{
        \int_{t_n}^{t} 
            \tfrac{2{v}_n-2Q_n I_n^h h
            -
            2\alpha {u}_nh
            }
            {2+\alpha h^2}
        {\rm{d}} s},
    \end{array}\right.
\end{equation}
for $t\in[t_n, t_{n+1}]$, $n \in \mathbb{N}_0$ and further denote
\begin{equation}\label{eq:def_Y_t}
    {Y}_t:= 
    \Big( 
        \begin{array}{c}
            \hat{v}_t
            \\
            \hat{u}_t 
        \end{array} 
    \Big),
    \
    t \geq 0.
\end{equation}
Before coming to the convergence analysis, we make an additional assumption as follows.
\begin{assumption}\label{ass:growth_Phi}
Suppose $\nabla \Phi \in C^2(\mathbb{R}^m; \mathbb{R}^m)$ and there exists some constant $C_{\Phi} \geq 0$ such that
$$
|
    \nabla \Phi(x)-\nabla \Phi(y)
|
\leq 
C_{\Phi}(1+|x|^2+|y|^2)|x-y|,
\
\forall x,y \in \mathbb{R}^m.
$$
\end{assumption}
Clearly,
Assumption \ref{ass:growth_Phi}
implies that 
for any $x=(x_1^*,x_2^*)^*,y=(y_1^*,y_2^*)^*$
with $x_1, x_2, y_1, y_2 \in \mathbb{R}^m$,
\begin{equation}\label{eq:|F(x)-F(y)|}
\begin{aligned}
|F(x)-F(y)| 
&=
\sqrt{
  |\kappa  (\nabla \Phi(x_2) - \nabla \Phi(y_2)) + \gamma (x_1 - y_1) |^2
  +
  | x_1 - y_1 |^2
}
\\
&\leq 
 \kappa| \nabla \Phi(x_2) - \nabla \Phi(y_2)| + (\gamma +1) |x_1 - y_1 |
\\
&\leq 
    \big(\kappa C_{\Phi}(1+|x_2|^2+|y_2|^2)+\gamma +1\big)|x-y|.
\end{aligned}
\end{equation}

\subsection{Order $\mathbf{\tfrac{1}{2}}$ strong convergence}
Since the exponential integrability property is only available for the numerical solution \eqref{eq:continu_ito_process_SSAV} on the mesh point $t_n$ (Theorem \ref{thm:exponential_integ_numer_solution}), instead of the entire interval $[0,T]$, we are not able to derive the order one strong convergence directly, relying on the perturbation estimates developed in \cite{dai2025perturbation}. 
Next we first show order $\mathbf{\tfrac{1}{2}}$ strong convergence of the numerical solution \eqref{eq:continu_ito_process_SSAV} on the mesh point.



\begin{proposition}\label{prop:half_and_one_convergence_rate}
    Let Assumptions \ref{ass:well_pose_rho_t} and \ref{ass:growth_Phi} hold and
assume
$\mathbb{E}\big[
        \exp(\delta(\tfrac{|v_0|^2}{2} 
        +\kappa \Phi(u_0)+C_H))
        \big]<+\infty$ and
    $\alpha {\delta}> e^{\lambda T} \kappa C_{\Phi} T \theta$ for some $\delta>0$.
     Let $\{X(t)\}_{t \in [0, T]}$ and the SSAV scheme $\{ Y_{t_n}\}_{n \in \{0,1,...,N_T\} }$ be defined by \eqref{eq:define_X_t_and_Y_n}, where $h\in (0,1]$.
    Then for any $\theta >0$, there exists a constant $C_T$ (exponentially depending on $T$) independent of $h$ such that 
    \begin{equation} \sup_{n\in\{0,1,...,N_T\}}
\mathbb{E}\big[ |X(t_{n})-{Y}_{t_n}|^\theta\big] \leq C_Th^{\frac{\theta}{2}}.
    \end{equation}
\end{proposition}
\textbf{Proof:} 
For $n \in \{0,1, ..., N_T\}$,
we first denote 
\begin{equation}\label{eq:denotes-e_n,J_n,K_n}
{e}_n:=X(t_n)-{Y}_{t_n},
\
{J}_n:=
\tfrac{-4Q_n I_n^h 
            -
            4\alpha {u}_n
            -
            2\alpha {v}_n h
            }
    {2+\alpha h^2}
,
\
\
\text{and}
\
{K}_n:=
\tfrac{2{v}_n-2Q_n I_n^h h
            -
            2\alpha {u}_nh
            }
    {2+\alpha h^2}.
\end{equation}
Then
the second and third equalities in \eqref{eq:solution_split_SAV}
will be simplified to 
$  v_n^{\vartriangle} = v_n + J_n h $
and
$ u_n^{\vartriangle} = u_n + K_n h $.
Inserting this into \eqref{eq:SSAV_OU_equation} and using \eqref{eq:langevin_sde_newform} we infer
\begin{equation}\label{eq:error-original-inLemma}
    \begin{aligned}
    &X(t_{n+1})-{Y}_{t_{n+1}}
    \\
    &\quad =
    \Big( 
                \begin{array}{c}
               e^{-\gamma h}v(t_n)
                \\
                u(t_n) 
                \end{array} 
    \Big)
    +
    \int_{t_n}^{t_{n+1}} 
            \Big( 
                \begin{array}{c}
                e^{-\gamma (t_{n+1}-s)}\big(-\kappa \nabla \Phi(u(s))\big) 
                \\
                v(s) 
                \end{array} 
             \Big)
    {\rm d}s
    +
    \int_{t_n}^{t_{n+1}} 
        \big( 
        \begin{array}{c}
         e^{-\gamma(t_{n+1}-s)}\Gamma
        \\
        \mathbf{0}
        \end{array} 
            \big) {\rm d}{W}_s
    \\
    &\qquad
    -   \Big( 
                \begin{array}{c}
               e^{-\gamma h}{v}_n
                \\
                {u}_n
                \end{array} 
    \Big)
    -
    \Big( 
        \begin{array}{c}
        e^{-\gamma h}{J}_n
        \\
        K_n
        \end{array}
    \Big) h
    -
       \int_{t_n}^{t_{n+1}} 
       \Big( 
        \begin{array}{c}
        e^{-\gamma(t_{n+1}-s)}\Gamma 
        \\
        \mathbf{0}
        \end{array} 
            \Big)
       {\rm{d}} W_{s}
\\
&\quad =
\underbrace{
 \Big( 
                \begin{array}{c}
               e^{-\gamma h}\big(v(t_n)-{v}_n\big)
                \\
                u(t_n)-{u}_n
                \end{array} 
    \Big)
}_{=:A_n}
+
\underbrace{
\int_{t_n}^{t_{n+1}} 
            \Big( 
                \begin{array}{c}
                 e^{-\gamma (t_{n+1}-s)}\kappa \big(-\nabla \Phi(u(t_n))\big)
                 -
                 e^{-\gamma h}{J}_n
                \\
                v(t_n) -K_n
                \end{array} 
             \Big)
    {\rm d}s
    }_{=:B_n}
\\
&\qquad +
\underbrace{
    \int_{t_n}^{t_{n+1}} 
            \Big( 
                \begin{array}{c}
                e^{-\gamma (t_{n+1}-s)}\kappa \big(-\nabla \Phi(u(s))+\nabla \Phi(u(t_n))\big)
                \\
                v(s) -v(t_n)
                \end{array} 
             \Big)
    {\rm d}s,
}_{=:{R}_n}
\end{aligned}
\end{equation}
where one can further split $B_n$ as 
\begin{equation}\label{eq:B_n=B_n1+B_n2+B_n3}
     \begin{aligned}
    B_n&=
            \int_{t_n}^{t_{n+1}}
           \Big( 
                \begin{array}{c}
                 (e^{-\gamma (t_{n+1}-s)}-1) \big(-\kappa \nabla 
                 \Phi(u(t_n))\big)
                 -
                (e^{-\gamma h}-1){J}_n
                \\
                0
                \end{array} 
             \Big)
             {\rm d}s
\\
    &\quad 
    +
            \Big( 
                \begin{array}{c}
                 \kappa \big(\nabla \Phi({u}_n)-
                 \nabla \Phi(u(t_n))\big) 
                \\
                v(t_n) -{v}_n
                \end{array} 
             \Big) h
    +
            \Big( 
                \begin{array}{c}
                -\kappa \nabla \Phi({u}_n) -{J}_n
                \\
                {v}_n -K_n
                \end{array} 
             \Big) h
\\
    &
    =: B_{n,1}+B_{n,2}+B_{n,3}.
    \end{aligned}
\end{equation}
With the aid of the Young inequality, one derives from \eqref{eq:error-original-inLemma} that
\begin{equation}
     \begin{aligned}
|{e}_{n+1}|^2 
    &=
    |A_n|^2+ |B_n|^2 + |{R}_n|^2+ 2 \langle  A_n, B_n \rangle+2 \langle A_n, {R}_n \rangle
    +2 \langle B_n, {R}_n \rangle
\\
    &\leq 
    |A_{n}|^2
    +2|B_n|^2+ 2|{R}_n|^2 
    + 2 \langle A_n, B_{n,2} \rangle
    +
    2h|A_{n}|^2
    +h^{-1}|B_{n,1}+B_{n,3}|^2
    +h^{-1}|{R}_n|^2 
\\
    &\leq 
    |{e}_{n}|^2
    +2h|{e}_{n}|^2 
    +2|B_n|^2+ 2|{R}_n|^2 
    + 2 | A_n| |B_{n,2}|
    +2h^{-1}|B_{n,1}|^2+2h^{-1}|B_{n,3}|^2
    +h^{-1}|{R}_n|^2 
\\
    &\leq 
    |{e}_{n}|^2
    +
    Ch|{e}_{n}|^2
     +
    2\kappa C_{\Phi} h(|u(t_n)|^2+|{u}_n|^2)|{e}_{n}|^2 
\\
    &\quad
     +2|B_n|^2
     + (2+h^{-1})|{R}_n|^2 
    +2h^{-1}|B_{n,1}|^2
    +2h^{-1}|B_{n,3}|^2 
    ,
    \end{aligned}
\end{equation} 
where the fact $|A_n| \leq |{e}_n|$ and the estimate \eqref{eq:|F(x)-F(y)|} 
were used in the second and the last inequality, respectively.
By iteration, we have
\begin{equation}
     \begin{aligned}
        |{e}_{n+1}|^2 
        &\leq  
        \sum_{k=0}^n 
        \big(
        Ch|{e}_{k}|^2
        +
        2 \kappa C_{\Phi} h(|u(t_k)|^2+|{u}_k|^2)|{e}_{k}|^2
        \big)
\\
    &\quad +
    \sum_{k=0}^n 
    \big(
    2|B_k|^2
     + (2+h^{-1})|{R}_k|^2 
    +2h^{-1}|B_{k,1}|^2
    +2h^{-1}|B_{k,3}|^2
    \big).
    \end{aligned}
\end{equation}
The discrete Gronwall inequality hence gives
\begin{equation}
     \begin{aligned}
        |{e}_{n+1}|^2 
    &\leq
    \Big(
     \sum_{k=0}^n
     \big(
        2|B_k|^2
        + (2+h^{-1})|{R}_k|^2 
        +2h^{-1}|B_{k,1}|^2
        +2h^{-1}|B_{k,3}|^2 
        \big)
        \Big)
     \\
        &\quad \cdot \exp\Big(
         2\kappa C_{\Phi}h
         \sum_{k=0}^n 
        (|u(t_k)|^2+|{u}_k|^2)
        +CT
        \Big).
    \end{aligned}
\end{equation}
Taking the norm $\| \cdot \|_{L^\frac{\theta}{2}(\Omega;\mathbb{R})}$ on both sides and using H\"{o}lder's inequality with $\tfrac{1}{\theta}=\tfrac{1}{p}+\tfrac{1}{q}$ for $p,q > \theta$, one arrives at
\begin{equation}
     \begin{aligned}
        \big\| {e}_{n+1} \big\|^2_{L^{\theta}(\Omega;\mathbb{R}^{2m})}
    &\leq
    \Big\|
     \sum_{k=0}^n
     \big(
       2|B_k|^2
        + (2+h^{-1})|{R}_k|^2 
        +2h^{-1}|B_{k,1}|^2
        +2h^{-1}|B_{k,3}|^2 
        \big)
        \Big\|_{L^{\frac{p}{2}}(\Omega;\mathbb{R})}
     \\
        &\quad \cdot \Big\|
        \exp\big(
         2\kappa C_{\Phi}h
         \sum_{k=0}^n 
        (|u(t_k)|^2+|{u}_k|^2)
        +CT
        \big)
        \Big\|_{L^{\frac{q}{2}}(\Omega;\mathbb{R})},
    \end{aligned}
\end{equation}
which in turn yields 
\begin{equation}\label{eq:e_n_in_Half_Order_Proof}
     \begin{aligned}
         \big\| {e}_{n+1} \big\|^2_{L^{\theta}(\Omega;\mathbb{R}^{2m})}
    &\leq
    C
    \sum_{k=0}^n
    \Big(
    \big\| B_k \big\|^2_{L^{p}(\Omega;\mathbb{R}^{2m})}
    + (1+h^{-1})\big\| {R}_k \big\|^2_{L^{p}(\Omega;\mathbb{R}^{2m})}
    + h^{-1} 
    \big\| B_{k,1} \big\|^2_{L^{p}(\Omega;\mathbb{R}^{2m})}
    \\
    &\quad \qquad \qquad
    +h^{-1} 
    \big\| B_{k,3} \big\|^2_{L^{p}(\Omega;\mathbb{R}^{2m})}
    \Big)
     \\
        &\quad \cdot 
        e^{CT}
        \Big\|
        \exp\big(
        \kappa
         C_{\Phi}h
         \sum_{k=0}^n 
        (|u(t_k)|^2+|{u}_k|^2)
        \Big\|^2_{L^{q}(\Omega;\mathbb{R})}.
    \end{aligned}
\end{equation}
Subsequently, it remains to bound terms $\| B_{k,1}\|_{L^{p}(\Omega;\mathbb{R}^{2m})}$,
$\| B_{k,2} \|_{L^{p}(\Omega;\mathbb{R}^{2m})}$
,
$\| B_{k,3} \|_{L^{p}(\Omega;\mathbb{R}^{2m})}$
and 
$\| R_{k}\|_{L^{p}(\Omega;\mathbb{R}^{2m})}$.
For the term $\| B_{k,1}\|_{L^{p}(\Omega;\mathbb{R}^{2m})}$, we first note
\begin{equation} \label{eq:estimate-Bk1}
\begin{aligned}
    &\big\| B_{k,1} \big\|_{L^{p}(\Omega;\mathbb{R}^{2m})}
\\
    &\quad \leq 
    \int_{t_k}^{t_{k+1}} 
    \Big(
    \kappa
     |e^{-\gamma (t_{n+1}-s)}-1| 
      \big\|
      \nabla \Phi(u(t_n))
      \big\|_{L^{p}(\Omega;\mathbb{R}^{m})}
     +
        |e^{-\gamma h}-1|
        \big\|
        {J}_n
        \big\|_{L^{p}(\Omega;\mathbb{R}^{m})}
    \Big)
    {\rm d}s.
\end{aligned}
\end{equation}
By Corollary \ref{coro:p_order_and_expone_moment} and Assumption \ref{ass:well_pose_rho_t}, one easily sees
\begin{equation}
\begin{aligned}
\big\|
      \nabla \Phi(u(t_n))
      \big\|_{L^{p}(\Omega;\mathbb{R}^{m})}
 {\ \textstyle \bigvee}
 \
\big\|
        {J}_n
        \big\|_{L^{p}(\Omega;\mathbb{R}^{m})}
        \leq C,
\end{aligned}
\end{equation}
which helps us derive from \eqref{eq:estimate-Bk1} that
\begin{equation}\label{eq:estimate_for_B_k1}
\begin{aligned}
 &\big\| B_{k,1} \big\|^2_{L^{p}(\Omega;\mathbb{R}^{2m})}
 \leq Ch^4.
\end{aligned}
\end{equation}
Owing to Corollary \ref{coro:p_order_and_expone_moment}, the term $\| B_{k,2} \|_{L^{p}(\Omega;\mathbb{R}^{2m})}$ can be similarly treated and bounded by
%
\begin{equation}\label{eq:estimate_for_B_k2_old}
\begin{aligned}
 &\big\| B_{k,2} \big\|^2_{L^{p}(\Omega;\mathbb{R}^{2m})}
 \leq C h^2.
\end{aligned}
\end{equation}
Now, let us concentrate on the estimate of 
$\| B_{k,3} \|_{L^{p}(\Omega;\mathbb{R}^{2m})}$. Recalling the definition of $K_n$, given in \eqref{eq:denotes-e_n,J_n,K_n}, we first show
\begin{equation}\label{eq:estimate_begin_B_K_3}
\begin{aligned}
    \big\| B_{k,3} \big\|^2_{L^{p}(\Omega;\mathbb{R}^{2m})}
    &\leq 
    2h^2
    \big(
    \big\|
        -\kappa \nabla \Phi({u}_k) -J_k
    \big\|^2_{L^p(\Omega;\mathbb{R}^{m})}
    +
    \big\|
      {v}_k -K_k 
    \big\|^2_{L^p(\Omega;\mathbb{R}^{m})}
    \big)
\\
    &\leq
    Ch^2
        \big\| -\kappa \nabla \Phi({u}_k) - J_k \big\|_{L^{p}(\Omega;\mathbb{R}^{m})}^2
       +
   C h^4.      
\end{aligned}
\end{equation}
Meanwhile, from \eqref{eq:I_n^h_definition} and \eqref{eq:express_rho_n_add_1} one infers ${\rho}_{n+1}+{\rho}_n=2I_n^h$, and further deduces
\begin{equation}
    J_k
    =
    \frac{ -2 Q_k ( {\rho}_{k+1} + {\rho}_k ) - 4 \alpha {u}_k - 2 \alpha {v}_k h }{2 + \alpha h^2}.
\end{equation}
This together with the notation \eqref{eq:Q_n_definition} shows
\begin{equation}
    \begin{aligned}
    &
       \big\| -\kappa \nabla \Phi({u}_k) - J_k \big\|_{L^{p}(\Omega;\mathbb{R}^{m})} 
\\
       &\quad =
       \big\| 
        2 Q_k  \sqrt{\kappa \Phi({u}_k)+C_H- \alpha |{u}_k|^{2}} + 2 \alpha {u}_k + J_k 
       \big\|_{L^{p}(\Omega;\mathbb{R}^{m})}
\\
       &\quad=
        \tfrac{1 }{  2 + \alpha h^2 } 
    \big\|      
    (4 +  2 \alpha h^2 ) Q_k \sqrt{\kappa \Phi({u}_k)+C_H- \alpha |{u}_k|^{2}}
        -
    2 Q_k ({\rho}_{k+1} + {\rho}_{k})
    +
    2\alpha^2  {u}_k h^2
    -
    2 \alpha {v}_k h
     \big\|_{L^{p}(\Omega;\mathbb{R}^{m})}
 \\
    &\quad=
        \tfrac{1 }{  2 + \alpha h^2 } 
    \Big\|      
    4 Q_k (\sqrt{\kappa \Phi({u}_k)+C_H- \alpha |{u}_k|^{2}}
    -
    {\rho}_k)
    +
        2 \alpha 
        Q_k \sqrt{\kappa \Phi({u}_k)+C_H- \alpha |{u}_k|^{2}}h^2
\\
    &\qquad \qquad \quad -
       Q_k \langle Q_k, v_k^{\vartriangle} + {v}_k \rangle h 
    +
    2\alpha^2 {u}_k h^2
    -
    2 \alpha {v}_k h
     \Big\|_{L^{p}(\Omega;\mathbb{R}^{m})}
\\
&\quad \leq
C
\big(
 \big\| \sqrt{\kappa \Phi({u}_k)+C_H- \alpha |{u}_k|^{2}} - {\rho}_k  \big\|_{L^{2p}(\Omega;\mathbb{R})} 
 +  h\big),
    \end{aligned}
\end{equation}
where the fact that
${\rho}_{k+1}={\rho}_k
    +\tfrac{1}{2}
   \langle Q_k, {v^{\vartriangle}_{k}+{v}_k}\rangle 
    h$  was utilized in the third equality and the last inequality stands due to Corollary   \ref{coro:p_order_and_expone_moment} and H\"{o}lder's inequality.
Next,
we claim that
\begin{equation}
\label{eq:claim_rho_k_Order_h}
    \big\| \sqrt{\kappa \Phi({u}_k)+C_H- \alpha |{u}_k|^{2}} - {\rho}_k  \big\|_{L^{2p}(\Omega;\mathbb{R})}
    \leq Ch,
\quad
\forall \ k \in \{0,1,...,N_T\}.
\end{equation}
By the Taylor expansion, we further use the notation \eqref{eq:Q_n_definition}, as well as the facts that
${\rho}_{k+1}={\rho}_k
    +\tfrac{1}{2}
   \langle Q_n, {v^{\vartriangle}_{k}+{v}_k}\rangle 
    h$ and ${u}_k -{u}_{k-1}=\tfrac{v^{\vartriangle}_{k-1}+{v}_{k-1}}{2}$,
to achieve
\begin{equation}
\begin{aligned}
&\sqrt{\kappa \Phi({u}_k)+C_H-\alpha |{u}_k|^2}
    - {\rho}_k
\\
&\quad =\sqrt{\kappa \Phi({u}_{k-1})+C_H-\alpha |{u}_{k-1}|^2}
    +
    \Big\langle 
    \int_0^1 
    \tfrac{\kappa \nabla \Phi({u}_{k-1}+r({u}_k-{u}_{k-1}))
            -2\alpha({u}_{k-1}+r({u}_k-{u}_{k-1}))}
    {2\sqrt{\kappa \Phi({u}_{k-1}+r({u}_k-{u}_{k-1}))+C_H
           -\alpha |{u}_{k-1}+r({u}_k -{u}_{k-1})|^2}}
    {\rm d}r,
    {u}_k-{u}_{k-1}
    \Big\rangle
\\
    &\qquad 
    -
    {\rho}_{k-1}
    -
    \tfrac{1}{2}
    \big\langle \tfrac{\kappa \nabla \Phi({u}_{k-1})-2\alpha {u}_{k-1}}{\sqrt{\kappa \Phi({u}_{k-1})+C_H- \alpha |{u}_{k-1}|^{2}}},
    \tfrac{v^{\vartriangle}_{k-1}+{v}_{k-1}}{2}\big\rangle 
    h
\\
    &\quad =\sqrt{\kappa \Phi({u}_{k-1})+C_H-\alpha |{u}_{k-1}|^2}-
    {\rho}_{k-1}
\\
    &\qquad 
    +
    \Big\langle
    \int_0^1  
    \Big(
    \tfrac{\kappa \nabla \Phi({u}_{k-1}+r({u}_k-{u}_{k-1}))
            -2\alpha({u}_{k-1}+r({u}_k-{u}_{k-1}))}
    {2\sqrt{\kappa \Phi({u}_{k-1}+r({u}_k-{u}_{k-1}))+C_H
           -\alpha |{u}_{k-1}+r({u}_k-{u}_{k-1})|^2}}
    -
    \tfrac{\kappa \nabla \Phi({u}_{k-1})-2\alpha {u}_{k-1}}{2\sqrt{\kappa \Phi({u}_{k-1})+C_H- \alpha |{u}_{k-1}|^{2}}}
    \Big)
    {\rm d}r
    ,
     \tfrac{v^{\vartriangle}_{k-1}+{v}_{k-1}}{2}
     \Big\rangle
     h,
\end{aligned}
\end{equation}
which ensures that
\begin{equation}
\begin{aligned}
&\big| \sqrt{\kappa \Phi({u}_k)+C_H-\alpha |{u}_k|^2}
    -{\rho}_k \big|
\\
&\quad \leq
h
\sum_{i=1}^{k}
\Big|
    \int_0^1  
    \Big(
    \tfrac{ \kappa \nabla \Phi({u}_{i-1}+r({u}_i-{u}_{i-1}))
            -2\alpha({u}_{i-1}+r({u}_i-{u}_{i-1}))}
    {2\sqrt{\kappa \Phi({u}_{i-1}+r({u}_i-{u}_{i-1}))+C_H
           -\alpha |{u}_{i-1}+r({u}_i-{u}_{i-1})|^2}}
    -
    \tfrac{\kappa \nabla \Phi({u}_{i-1})-2\alpha {u}_{i-1}}{2\sqrt{\kappa \Phi({u}_{i-1})+C_H- \alpha |{u}_{i-1}|^{2}}}
    \Big)
    {\rm d}r
\Big|
\cdot
\Big|
     \tfrac{v^{\vartriangle}_{n-1}+{v}_{n-1}}{2}
\Big|.
\end{aligned}
\end{equation}
Again, using the Taylor expansion and the fact that $u_{n+1} = u_n + \tfrac{{v}^{\vartriangle}_{n}+{v}_n}{2}h, n\in \{0,1,...,N_T-1\}$, together with Assumption \ref{ass:growth_Phi} and 
Corollary  \ref{coro:p_order_and_expone_moment},
gives the claim \eqref{eq:claim_rho_k_Order_h}.
This in turn gives
  \begin{equation}\label{eq:-nablaPHI-J-inproof}
    \big\| -\kappa \nabla \Phi({u}_k) - J_k \big\|_{L^{p}(\Omega;\mathbb{R}^{m})}  \leq
       C h,
  \end{equation}
and it is thus concluded that
\begin{equation}\label{eq:estimate_for_B_k3}
\begin{aligned}
    \big\| B_{k,3} \big\|^2_{L^{p}(\Omega;\mathbb{R}^{2m})}
    &\leq 
   C h^4.      
\end{aligned}
\end{equation}
When it comes to the term
$\| {R}_{k}\|_{L^{p}(\Omega;\mathbb{R}^{2m})}$,
by noting that the process $\{X(t)\}_{t\in[0,T]}$ has a mean-square H\"{o}lder continuity of order $\tfrac{1}{2}$, one easily obtains
\begin{equation}\label{eq:estimate_for_R_k_SSAV1}
    \begin{aligned}
      \big\| {R}_k \big\|^2_{L^{p}(\Omega;\mathbb{R}^{2m})}
      &\leq
      C h^3.
    \end{aligned}
\end{equation}
Finally, we deal with the remaining term in \eqref{eq:e_n_in_Half_Order_Proof}.
In view of Jensen's inequality (see, e.g., \cite[Lemma 2.22]{cox2024local}) and Corollary \ref{coro:p_order_and_expone_moment}, one can show
\begin{equation}\label{eq:estimate_for_expon_term} 
     \begin{aligned}
        &\Big\|
        \exp\big(
         \kappa C_{\Phi}h
        \sum_{k=0}^{N-1}
        (|u(t_k)|^2+|{u}_k|^2)\big)
        \Big\|_{L^{q}(\Omega;\mathbb{R})}
    \\
        &\quad \leq 
        \tfrac{1}{T} h{}
        \sum_{k=0}^{N-1}
         \big\|
        \exp\big(
         \kappa C_{\Phi}T
        (|u(t_k)|^2+|{u}_k|^2)\big)
        \big\|_{L^{q}(\Omega;\mathbb{R})}
    \\
    &\quad \leq 
    \sup_{k=0,...,N-1}  
    \big\|
        \exp\big(
        \kappa C_{\Phi}T
        (|u(t_k)|^2+|{u}_k|^2)\big)
    \big\|_{L^{q}(\Omega;\mathbb{R})}
    \\
    &\quad \leq C,
    \end{aligned}
\end{equation}
where the assumption
$ \alpha \delta \geq e^{\lambda T}\kappa C_{\Phi}Tq$ was used.
Gathering \eqref{eq:estimate_for_B_k1},
\eqref{eq:estimate_for_B_k2_old} and \eqref{eq:estimate_for_B_k3}-\eqref{eq:estimate_for_expon_term},
one derives from \eqref{eq:e_n_in_Half_Order_Proof} that
\begin{equation}
     \begin{aligned}
       \sup_{n\in\{0,1,...,N_T\}} \big\| {e}_{n} \big\|_{L^{\theta}(\Omega;\mathbb{R}^{2m})}
    &\leq  C_Th^{\frac{1}{2}},
    \end{aligned}
\end{equation}
which completes the proof.
\qed

As a direct consequence of Corollary   \ref{coro:p_order_and_expone_moment} and Proposition \ref{prop:half_and_one_convergence_rate}, the following corollary provides a half order convergence result for the continuous version of the SSAV method.

\begin{corollary}\label{coro:half_convergence_rate}
  Let the assumptions of Proposition \ref{prop:half_and_one_convergence_rate} be fulfilled.
    Let $ \{X(t)\}_{t \in [0,T]}$ be defined by \eqref{eq:langevin_sde} and 
 $\{ Y_t\}_{t \in [0,T]}$ be given by \eqref{eq:def_Y_t}.
For any $p>0$, there exists a constant $C_T$ (exponentially depending on $T$) independent of $h$, such that 
    \begin{equation}
\sup_{t\in[0,T]}
\mathbb{E}[| X(t) - Y_t |^p] 
\leq 
C_T h^{\frac{p}{2}}.
    \end{equation}
\end{corollary}

\textbf{Proof:} In view of \eqref{eq:continu_ito_process_SSAV} and by Corollary  \ref{coro:p_order_and_expone_moment}, one easily deduces 
\begin{equation}
    \mathbb{E}[|Y_{\lfloor t \rfloor_h}-Y_{t}|^p]
    \leq
    C h^{\frac{p}{2}}, 
    \quad
    \forall t \in [0,T],
\end{equation}
which together with 
Proposition
\ref{prop:half_and_one_convergence_rate}
implies
\begin{equation}
\begin{aligned}
&
\sup_{t\in[0,T]}
\mathbb{E}[|X(t)-Y_t|^p]
\\
&\quad \leq 
C_{p}
\big(
    \sup_{t\in[0,T]}
    \mathbb{E}[|X(t)-X(\lfloor t \rfloor_h)|^p]
    +
    \sup_{t\in[0,T]}
    \mathbb{E}[|X(\lfloor t \rfloor_h)-Y_{\lfloor t \rfloor_h}|^p]
    +
    \sup_{t\in[0,T]}
    \mathbb{E}[|Y_{\lfloor t \rfloor_h}-Y_{t}|^p]
\big)
\\
&\quad \leq Ch^{\frac{p}{2}},
\end{aligned}
\end{equation}
which gives the assertion.\qed

\subsection{Order one strong convergence}
In this subsection, we aim to derive the desired order one convergence of the proposed scheme, based on the previously obtained order $\tfrac{1}{2}$ convergence. For convenience, in the remainder of this paper we further denote
\begin{equation}
a(t) :=\Big( 
        \begin{array}{c}
            \tfrac{-4Q_n I_n^h 
            -
            4\alpha u_n
            -
            2\alpha v_n h
            }
    {2+\alpha h^2}
    -\gamma \hat{v}_t 
            \\
     \tfrac{2v_n-2Q_n I_n^h h
            -
            2\alpha u_n h
            }
            {2+\alpha h^2}
        \end{array} 
    \Big),
\ 
b :=\Big( 
\begin{array}{c}
    \Gamma 
     \\
    \mathbf{0}
    \end{array} 
            \Big)=G,
\end{equation}
and
\begin{equation}
{f}(x):=
    \Big( 
        \begin{array}{c}
            -\kappa \nabla \Phi(x_2)
            \\
            x_1 
        \end{array} 
    \Big),\
    x_1, x_2\in \mathbb{R}^m.
\end{equation}
Following this notation, $\{Y_{t}\}_{t\geq 0}$, given by \eqref{eq:continu_ito_process_SSAV}-\eqref{eq:def_Y_t}, can be rewritten as
\begin{equation}
{\rm d} Y_t=
a(t)\, {\rm d}t + b \, {\rm d}W_t, \ t> 0.
\end{equation}
Equipped with Corollary \ref{coro:half_convergence_rate}, we are now well-prepared to show the optimal strong convergence rate for the SSAV method in a finite time horizon.
To this end, we introduce a novel auxiliary process $\{\chi_t\}_{t \in [0,T]}$, given by \eqref{eq:chi-t}, different from those in \cite{Hutzenthaler2020,dai2025perturbation}.

\begin{theorem}[Optimal strong convergence rate]
\label{thm:optimal_converge_rate}
 Let Assumptions \ref{ass:well_pose_rho_t},
 \ref{ass:growth_Phi} hold. For $\delta, \theta>0$,
assume that 
$\mathbb{E}\big[
        \exp(\delta(\tfrac{|v_0|^2}{2} +\kappa \Phi(u_0)+C_H))
        \big]<+\infty$ and 
    $\alpha {\delta}> e^{\lambda T}\kappa C_{\Phi} T \theta  $.
     Let $ \{X(t)\}_{
    t \in [0,T]} $ and the SSAV scheme $ \{ Y_{t_n}\}_{n \in \{0,1,...,N_T\} }$ be defined by \eqref{eq:define_X_t_and_Y_n}, where $h\in(0,1]$. 
    Then
\begin{equation}
\mathbb{E}\Big[\sup_{n\in\{0,1,...,N_T\}}|X(t_{n})-Y_{t_n}|^{\theta}\Big] \leq C_Th^{\theta},
\end{equation}
where $C_T$ is some constant independent of $h$, but exponentially depending on $T$.
\end{theorem}
\textbf{Proof:}
For brevity, we introduce
\begin{equation} \label{eq:chi-t}
 \chi_t := 
 2\big(\kappa C_{\Phi}(1+|u({t_n})|^2+|u_n|^2)+\gamma +1\big)
 , \  t\in[t_n, t_{n+1}),
 \quad
 n \in \{0,1,...,N_T-1\}.
\end{equation}
Using the It\^{o} formula and \eqref{eq:|F(x)-F(y)|} then gives
\begin{equation}
\begin{aligned}
        \tfrac{|X(t) - Y_t|^2}
            {\exp (\int_0^t \chi_s {\rm d}s)} 
        &=
        |X(0) - Y_0|^2+ 
        \int_0^t 
            \tfrac {2
                \langle 
                    X(s) - Y_s,\left(G-b\right){\rm d} W_s 
                \rangle  }{\exp (\int_0^s \chi_r {\rm d}r)}
        \\
        &\quad +\int_0^t 
                    \tfrac{2 \langle X(s) - Y_s,F(X(s))-a(s)\rangle  
                        -
                        |X(s) - Y_s|^2\chi_s}
                    {\exp (\int_0^s \chi_r {\rm d}r)}
                {\rm d}s
            +\int_0^t 
        \tfrac{\| G-b \|^2}
        {\exp (\int_0^s \chi_r {\rm d}r)} {\rm d}s
        \\
        &=
        \int_0^t 
                    \tfrac{2 \langle X(s) - Y_s,F(X(s))-F(Y_s)\rangle  
                        -
                        |X(s) - Y_s|^2\chi_s}
                    {\exp (\int_0^s \chi_r {\rm d}r)}
                {\rm d}s
        +
        \int_0^t 
                    \tfrac{2 \langle X(s) - Y_s, F(Y_s)-a(s)\rangle } 
                    {\exp (\int_0^s \chi_r {\rm d}r)}
                {\rm d}s
\\
    &\leq 
    \int_0^t 
                    \tfrac{2\kappa C_{\Phi} |X(s)-Y_s|^2
                        \left( 
                            |u(s)|^2
                            -
                            |u( \lfloor s \rfloor_h)|^2
                            +
                            |\hat{u}_s|^2
                            -
                            |\hat{u}_{\lfloor s \rfloor_h}|^2
                          \right)}
                    {\exp (\int_0^s \chi_r {\rm d}r)}
                {\rm d}s
        +
        \int_0^t 
                    \tfrac{2 \langle X(s) - Y_s, F(Y_s)-a(s)\rangle } 
                    {\exp (\int_0^s \chi_r {\rm d}r)}
                {\rm d}s.
        \end{aligned}
      \end{equation}
Owing to \eqref{eq:langevin_sde} and \eqref{eq:continu_ito_process_SSAV}, we clearly have
$\| u(t)- u(\lfloor t \rfloor_h) \|_{L^{\frac{p}{2}}(\Omega; \mathbb{R}^{m})} + \| \hat{u}_t - \hat{u}_{\lfloor t \rfloor_h} \|_{L^{\frac{p}{2}}(\Omega; \mathbb{R}^{m})} \leq C h$, which together with Corollary \ref{coro:half_convergence_rate} 
yields, for any $ \vartheta \in [0,T]$,
\begin{equation}\label{eq:order_one_estima_origin}
\begin{aligned}
    &\Big\| \sup_{t\in[0,\vartheta]}
            \tfrac{|X(t) - Y_t|^2}
            {\exp (\int_0^t \chi_s {\rm d}s)}
    \Big\|_{L^{\frac{p}{2}}(\Omega; \mathbb{R})} 
\\
    &\quad \leq 
    \int_0^\vartheta
    \big\| 
        \tfrac{2 \kappa C_{\Phi} |X(t)-Y_t|^2
                \left(
                |u(t) + u(\lfloor t \rfloor_h)|
                \cdot
                |u(t) - u(\lfloor t \rfloor_h)|
                    +
                    |\hat{u}_t + \hat{u}_{\lfloor t \rfloor_h}|
                    \cdot
                    |\hat{u}_t - \hat{u}_{\lfloor t \rfloor_h}|
                \right)}
        {\exp (\int_0^t \chi_s {\rm d}s)}
    \big\|_{L^{\frac{p}{2}}(\Omega; \mathbb{R})} 
    {\rm d}t
\\
    &\qquad +
    \Big\|
        \sup_{t\in[0,\vartheta]}\int_0^t 
            \tfrac{2 \langle X(s) - Y_s, F(Y_s)-a(s)\rangle } 
            {\exp (\int_0^s \chi_r {\rm d}r)}
        {\rm d}s
        \Big\|_{L^{\frac{p}{2}}(\Omega; \mathbb{R})}
\\
&\quad \leq 
\Big\|
        \sup_{t\in[0,\vartheta]}\int_0^t 
            \tfrac{2 \langle X(s) - Y_s, F(Y_s)-a(s)\rangle } 
            {\exp (\int_0^s \chi_r {\rm d}r)}
        {\rm d}s
        \Big\|_{L^{\frac{p}{2}}(\Omega; \mathbb{R})}
+Ch^2.
        \end{aligned}
      \end{equation}
Following the notation \eqref{eq:B_n=B_n1+B_n2+B_n3},
one can write,
for $s \in [t_n, t_{n+1}), n \in \{0,1,...,N_T-1\}$,
\begin{equation}\label{eq:decompo_f_Y_s-a_s}
\begin{aligned}
     F(Y_s)-a(s)
     &=
    \Big( 
        \begin{array}{c}
            -\kappa \nabla \Phi(\hat{u}_s) - \gamma \hat{v}_s
            \\
            \hat{v}_s
        \end{array} 
    \Big)
    -
    \Big( 
        \begin{array}{c}
            \tfrac{-4Q_n I_n^h 
            -
            4\alpha u_n
            -
            2\alpha v_n h
            }
    {2+\alpha h^2}
    -\gamma \hat{v}_s 
            \\
     \tfrac{2v_n-2Q_n I_n^h h
            -
            2\alpha u_nh
            }
            {2+\alpha h^2}
        \end{array} 
    \Big)
\\
&={{f}(Y_s)-{f}(Y_{t_n})}
    +
    h^{-1} B_{n,3}.
        \end{aligned}
\end{equation}
Further, we apply It\^{o}'s formula and acquire
\begin{equation}\label{eq:f(Y_s)-f(Y_n)-ItoFormula}
\begin{aligned}
f(Y_s)-f(Y_{t_n})
=
\int_{t_n}^s 
\big(
f'(Y_r)a(r)
    +
\tfrac{1}{2}\sum_{j=1}^m    
    f''(Y_r)(be_j,be_j)
\big) 
 {\rm d}r
+
\int_{t_n}^s 
f'(Y_r)b 
{\rm d}W_r,
        \end{aligned}
\end{equation}
where $\{ e_j \}_{ j \in \{1,2,...,m\} }$ is the orthonormal basis of $\mathbb{R}^m$. 
By the H\"{o}lder inequality and the Young inequality,
one can deduce
\begin{equation}\label{eq:X(s)-Ys,F(Ys)-a(s)_in_Thm}
\begin{aligned}
    &
    \Big\|
        \sup_{t\in[0,\vartheta]}\int_0^t 
            \tfrac{2 \langle X(s) - Y_s, F(Y_s)-a(s)\rangle } 
            {\exp (\int_0^s \chi_r {\rm d}r)}
        {\rm d}s
    \Big\|_{L^{\frac{p}{2}}(\Omega; \mathbb{R})}
\\
&\quad\leq 
\Big\|
        \sup_{t\in[0,\vartheta]}\int_0^t 
            \tfrac{2 \left\langle X(s) - Y_s, 
            \int_{\lfloor s \rfloor_h}^s
            \left(
                    f'(Y_r)a(r)
                    +
                    \frac{1}{2}\sum_{j=1}^m    
                    f''(Y_r)(be_j,be_j)
            \right) 
                    {\rm d}r
            \right\rangle } 
            {\exp (\int_0^s \chi_r {\rm d}r)}
        {\rm d}s
    \Big\|_{L^{\frac{p}{2}}(\Omega; \mathbb{R})}
\\
&\qquad 
+
    \Big\|
        \sup_{t\in[0,\vartheta]}\int_0^t 
            \tfrac{2 \left\langle X(s) - Y_s, 
                \int_{\lfloor s \rfloor_h}^s 
                f'(Y_r)b 
                    {\rm d}W_r
            \right\rangle } 
            {\exp (\int_0^s \chi_r {\rm d}r)}
        {\rm d}s
    \Big\|_{L^{\frac{p}{2}}(\Omega; \mathbb{R})}
+   
    \Big\|
        \sup_{t\in[0,\vartheta]}\int_0^t 
            \tfrac{2 \left\langle X(s) - Y_s, B_{\lfloor s \rfloor_h/h,3} 
            \right\rangle } 
            {\exp (\int_0^s \chi_r {\rm d}r)}
           \cdot h^{-1}
        {\rm d}s
    \Big\|_{L^{\frac{p}{2}}(\Omega; \mathbb{R})}
\\
&\quad \leq 
C
    \int_0^\vartheta
        \Big\|
                \tfrac{|X(s) - Y_s|} 
                {\exp (\int_0^s \frac12\chi_r {\rm d}r)}
        \Big\|^2_{L^{p}(\Omega; \mathbb{R})}
        {\rm d}s
    +
        \Big\|
        \sup_{t\in[0,\vartheta]}\int_0^t 
            \tfrac{2 \left\langle X(s) - Y_s, 
                \int_{\lfloor s \rfloor_h}^s 
                f'(Y_r)b 
                    {\rm d}W_r
            \right\rangle } 
            {\exp (\int_0^s \chi_r {\rm d}r)}
        {\rm d}s
        \Big\|_{L^{\frac{p}{2}}(\Omega; \mathbb{R})}
    +Ch^2,
        \end{aligned}
      \end{equation}
where Corollary   \ref{coro:p_order_and_expone_moment} and \eqref{eq:estimate_for_B_k3} were also used.
Thanks to It\^{o}'s formula, one notes that
\begin{equation}
        \begin{aligned} 
            \tfrac{X(s) - Y_s}{\exp (\int_0^s {\chi_r} {\rm d}r)}
            &=
            \tfrac{X({\lfloor s \rfloor_h}) - Y_{\lfloor s \rfloor_h}}
            {\exp (\int_0^{\lfloor s \rfloor_h} {\chi_r} {\rm d}r)}
            +
             \int_{\lfloor s \rfloor_h}^s \tfrac{F(X(r)) - a(r)}{\exp (\int_0^r {\chi_{\iota}} {\rm d}\iota)}{\rm d}r
            +
            \int_{\lfloor s \rfloor_h}^s \tfrac{(X(r)-Y_r)(-\chi_r)}{\exp (\int_0^r {{\chi_{\iota}}} {\rm{d}}\iota)}{\rm d}r,
        \end{aligned}
    \end{equation}
which then admits the following decomposition:
\begin{equation}
\begin{aligned}
&\Big\|
    \sup_{t\in[0,\vartheta]}\int_0^t 
            \tfrac{2 \left\langle X(s) - Y_s, 
                \int_{\lfloor s \rfloor_h}^s 
                f'(Y_r)b 
                    {\rm d}W_r
            \right\rangle } 
            {\exp (\int_0^s \chi_r {\rm d}r)}
        {\rm d}s
        \Big\|_{L^{\frac{p}{2}}(\Omega; \mathbb{R})} 
\\
&\quad
    \leq  
    \Big\|
    2\sup_{t\in[0,\vartheta]}\int_0^t 
    \Big\langle
            \tfrac{X(\lfloor s \rfloor_h) - Y_{\lfloor s \rfloor_h}}
            {\exp (\int_0^{\lfloor s \rfloor_h} {\chi_r} {\rm d}r)}
            ,
                \int_{\lfloor s \rfloor_h}^s 
                f'(Y_r)b 
                    {\rm d}W_r
        \Big\rangle
        {\rm d}s
        \Big\|_{L^{\frac{p}{2}}(\Omega; \mathbb{R})}
\\
    &\qquad 
    +
    \Big\|
    2\sup_{t\in[0,\vartheta]}\int_0^t 
    \Big\langle
             \int_{\lfloor s \rfloor_h}^s \tfrac{F(X(r)) - a(r)}{\exp (\int_0^r {\chi_{\iota}} {\rm d}\iota)}{\rm d}r
            ,
                \int_{\lfloor s \rfloor_h}^s 
                f'(Y_r)b 
                    {\rm d}W_r
        \Big\rangle
        {\rm d}s
        \Big\|_{L^{\frac{p}{2}}(\Omega; \mathbb{R})}
\\
    &\qquad 
    +
    \Big\|
    2\sup_{t\in[0,\vartheta]}\int_0^t 
    \Big\langle
             \int_{\lfloor s \rfloor_h}^s \tfrac{(X(r)-Y_r)(-\chi_r)}{\exp (\int_0^r {{\chi_{\iota}}} {\rm{d}}\iota)}{\rm d}r
            ,
                \int_{\lfloor s \rfloor_h}^s 
                f'(Y_r)b 
                    {\rm d}W_r
        \Big\rangle
        {\rm d}s
        \Big\|_{L^{\frac{p}{2}}(\Omega; \mathbb{R})}
\\
&\quad =: S_1+S_2+S_3.
        \end{aligned}
\end{equation}
To proceed, we split the term $S_1$ further:
\begin{equation}
\begin{aligned}
\| S_1 \|_{L^{\frac{p}{2}}(\Omega; \mathbb{R})}
&=
\Big\|
    2\sup_{t\in[0,\vartheta]}\int_0^t 
    \Big\langle
            \tfrac{X(\lfloor s \rfloor_h) - Y_{\lfloor s \rfloor_h}}
            {\exp (\int_0^{\lfloor s \rfloor_h} {\chi_r} {\rm d}r)}
            ,
                \int_{\lfloor s \rfloor_h}^s 
                f'(Y_r)b 
                    {\rm d}W_r
        \Big\rangle
        {\rm d}s
        \Big\|_{L^{\frac{p}{2}}(\Omega; \mathbb{R})}
\\
&\leq
    \Big\| 
    \sup_{t\in[0,\vartheta]} 
        \Big|
            \sum\limits_{k=0}^{\lfloor t \rfloor_h/h-1}\int_{t_k}^{t_{k+1}} 
                {2\Big\langle 
                    \tfrac{X(\lfloor s \rfloor_h) - Y_{\lfloor s \rfloor_h}}
                    {\exp (\int_0^{\lfloor s \rfloor_h} {\chi_r} {\rm d}r)},\int_{\lfloor s \rfloor_h}^s 
                        f'(Y_r)b{\rm d}W_r
                    \Big\rangle }{\rm{d}}s
        \Big|
    \Big\|_{L^{\frac{p}{2}}(\Omega; \mathbb{R})}
\\
    &\quad 
    +
    \Big\|\sup_{t\in[0,\vartheta]}
        \int_{\lfloor t \rfloor_h}^t 
        {2\Big\langle 
            \tfrac{X(\lfloor s \rfloor_h) - Y_{\lfloor s \rfloor_h}}
            {\exp (\int_0^{\lfloor s \rfloor_h} {\chi_r} {\rm d}r)},\int_{\lfloor s \rfloor_h}^s 
            f'(Y_r)b{\rm{d}}W_r 
         \Big\rangle }{\rm d }s
    \Big\|_{L^{\frac{p}{2}}(\Omega; \mathbb{R})}
\\
&=:S_{1,1}+S_{1,2},
        \end{aligned}
\end{equation}
where, by the It\^o isometry and Lemma \ref{prop:half_and_one_convergence_rate},
the term $S_{1,2}$ can be immediately estimated by
\begin{equation}\label{eq:order_one_estima_S_12}
    \begin{aligned}
        S_{1,2}\leq Ch^2.
    \end{aligned}
\end{equation}
In contrast, the term $S_{1,1}$ should be treated more carefully. 
Noting that 
    \[
        \zeta_n
        : =
        \sum\limits_{k=0}^{n-1}\int_{t_k}^{t_{k+1}} 
        {2\Big\langle \tfrac{X(\lfloor s \rfloor_h) - Y_{\lfloor s \rfloor_h}}
        {\exp (\int_0^{\lfloor s \rfloor_h} {\chi_r} {\rm d}r)}
        ,
            \int_{\lfloor s \rfloor_h}^s 
            f'(Y_r)b{\rm d}W_r 
        \Big\rangle }{\rm d}s,\ n\in \mathbb{N}
    \]
    is a discrete martingale,
    the Doob discrete martingale inequality, the Burkholder-Davis-Gundy inequality (see, e.g., \cite[Lemma 4.1]{hutzenthaler2011convergence}) and H{\"o}lder's inequality then give
\begin{equation}\label{eq:order_one_estima_S_11}
    \begin{aligned}
    S_{1,1} 
    &\leq 
    C_p \Big\| 
        \sum\limits_{k=0}^{\lfloor \vartheta \rfloor_h/h-1}\int_{t_k}^{t_{k+1}} 
        {2\Big\langle 
            \tfrac{X(\lfloor s \rfloor_h) - Y_{\lfloor s \rfloor_h}}
            {\exp (\int_0^{\lfloor s \rfloor_h} {\chi_r} {\rm d}r)},\int_{\lfloor s \rfloor_h}^s 
            f'(Y_r)b{\rm{d}}W_r 
         \Big\rangle }{\rm d }s
    \Big\|_{L^{\frac{p}{2}}(\Omega; \mathbb{R})}
\\
    &\leq C_p \bigg(\sum\limits_{k=0}^{\lfloor \vartheta \rfloor_h/h-1} 
    \Big\| 
    \int_{t_k}^{t_{k+1}} 
        {2\Big\langle 
            \tfrac{X(\lfloor s \rfloor_h) - Y_{\lfloor s \rfloor_h}}
            {\exp (\int_0^{\lfloor s \rfloor_h} {\chi_r} {\rm d}r)},\int_{\lfloor s \rfloor_h}^s 
            f'(Y_r)b{\rm{d}}W_r 
         \Big\rangle }{\rm d }s
    \Big\|^2_{{L^{\frac{p}{2}}(\Omega;\mathbb{R})}}\bigg)^{\frac{1}{2}}
\\
    &\leq C_p 
    \bigg(h \int_{0}^{\vartheta} 
    \Big\|
        \tfrac{|X(\lfloor s \rfloor_h)-Y_{\lfloor s \rfloor_h}|}
        {\exp (\int_0^{\lfloor s \rfloor_h} 
        { \frac12 }
        {{\chi_r}} {\rm d}r)}
    \Big\|^2_{{L^{p}(\Omega;\mathbb{R})}} 
    \Big\|
            \int_{\lfloor s \rfloor_h}^s  f'(Y_r)b{\rm{d}}W_r
    \Big\|^2_{{L^{p}(\Omega;\mathbb{R}^{2m})}}{\rm{d}}s
    \bigg)^{\frac{1}{2}}
\\
    &\leq 
        C_p
        \sup_{t\in [0,\vartheta]} 
        \Big\|
        \tfrac{|X(t)-Y_t|}{\exp (\int_0^t { \frac12 } {\chi_r} {\rm d}r)}
        \Big\|_{L^p(\Omega;\mathbb{R})}
         \bigg(
                h \int_{0}^{\vartheta}
                \Big\|
                    \int_{\lfloor s \rfloor_h}^s  
                    f'(Y_r)b{\rm d}
                    W_r
                \Big\|^2_{L^p(\Omega;\mathbb{R}^{2m})}{\rm d}s
            \bigg)^{\frac{1}{2}}
\\
            &\leq 
            \tfrac{1}{8}
            \sup_{t\in [0,\vartheta]} 
            \Big\|
            \tfrac{|X(t)-Y_t|}{\exp (\int_0^t { \frac12 } {\chi_r} {\rm d}r)}
            \Big\|^2_{L^p(\Omega;\mathbb{R})}
            +Ch^2.
    \end{aligned}
\end{equation}
Next, we turn to the term $S_2$, which, by virtue of the equations \eqref{eq:|F(x)-F(y)|}, \eqref{eq:f(Y_s)-f(Y_n)-ItoFormula}, and the estimate \eqref{eq:estimate_for_B_k3}, is estimated as follows:
\begin{equation}\label{eq:order_one_estima_S_2}
\begin{aligned}
S_2 &\leq 
\Big\|
    2\sup_{t\in[0,\vartheta]}\int_0^t 
    \Big\langle
             \int_{\lfloor s \rfloor_h}^s \tfrac{F(X(r)) - F(Y_r)}{\exp (\int_0^r {\chi_{\iota}} {\rm d}\iota)}{\rm d}r
            ,
                \int_{t_n}^s 
                f'(Y_r)b 
                    {\rm d}W_r
        \Big\rangle
        {\rm d}s
        \Big\|_{L^{\frac{p}{2}}(\Omega; \mathbb{R})}
\\
&\quad +
\Big\|
    2\sup_{t\in[0,\vartheta]}\int_0^t 
    \Big\langle
             \int_{\lfloor s \rfloor_h}^s \tfrac{ F(Y_r)-a(r)}{\exp (\int_0^r {\chi_{\iota}} {\rm d}\iota)}{\rm d}r
            ,
                \int_{t_n}^s 
                f'(Y_r)b 
                    {\rm d}W_r
        \Big\rangle
        {\rm d}s
        \Big\|_{L^{\frac{p}{2}}(\Omega; \mathbb{R})}
\\
&\leq
    C\int_0^\vartheta
        \Big(
             \int_{\lfloor s \rfloor_h}^s 
       \Big\|
             \tfrac{
             ( 1 + |u(r)|^2 + |\hat{u}_r|^2 )
             | X(r) - Y_r |
             }{\exp (\int_0^r {\chi_{\iota}} {\rm d}\iota)}
        \Big\|_{L^{p}(\Omega; \mathbb{R})}
             {\rm d}r
        \Big)
    \cdot
    \Big\|
           \int_{t_n}^s 
          f'(Y_r)b 
                    {\rm d}W_r
        \Big\|_{L^{p}(\Omega; \mathbb{R}^{2m})}
        {\rm d}s
\\
&\quad +
    2
    \int_0^\vartheta
       \Big(
       \int_{\lfloor s \rfloor_h}^s 
       \Big\|
       \tfrac{ 
       |f(Y_r) - f(Y_{\lfloor r \rfloor_h})|
       +
       {h^{-1}}| B_{\lfloor r \rfloor_h/h,3} | }
       {\exp (\int_0^r {\chi_{\iota}} {\rm d}\iota)}
       \Big\|_{L^{p}(\Omega; \mathbb{R})}
       {\rm d}r
       \Big)
        \cdot
        \Big\|
                \int_{t_n}^s 
                f'(Y_r)b 
                    {\rm d}W_r
        \Big\|_{L^{p}(\Omega; \mathbb{R}^{2m})}
            {\rm d}s
\\
&\leq 
Ch^2.
\end{aligned}
\end{equation}
In a similar manner, we acquire
\begin{equation}\label{eq:order_one_estima_S_3}
\begin{aligned}
S_3 &\leq Ch^2.
\end{aligned}
\end{equation}
Taking \eqref{eq:order_one_estima_S_11}, \eqref{eq:order_one_estima_S_12}, \eqref{eq:order_one_estima_S_2} and \eqref{eq:order_one_estima_S_3} into account, one obtains 
\begin{equation}
\Big\|
    \sup_{t\in[0,\vartheta]}\int_0^t 
            \tfrac{2 \left\langle X(s) - Y_s, 
                \int_{t_n}^s 
                f'(Y_r)b 
                    {\rm d}W_r
            \right\rangle } 
            {\exp (\int_0^s \chi_r {\rm d}r)}
        {\rm d}s
        \Big\|_{L^{\frac{p}{2}}(\Omega; \mathbb{R})} 
\leq 
\tfrac{1}{8}
            \sup_{t\in [0,\vartheta]} 
            \Big\|
            \tfrac{|X(t)-Y_t|}{\exp (\int_0^t \frac{1}{2}{\chi_r} {\rm d}r)}
            \Big\|^2_{L^p(\Omega;\mathbb{R})}
            +Ch^2.
\end{equation}
Inserting the above estimate into \eqref{eq:X(s)-Ys,F(Ys)-a(s)_in_Thm},
we arrive at
\begin{equation}
\begin{aligned}
    &
    \Big\|
        \sup_{t\in[0,\vartheta]}\int_0^t 
            \tfrac{2 \langle X(s) - Y_s, F(Y_s)-a(s)\rangle } 
            {\exp (\int_0^s \chi_r {\rm d}r)}
        {\rm d}s
    \Big\|_{L^{\frac{p}{2}}(\Omega; \mathbb{R})}
\\
&\quad \leq 
C
    \int_0^\vartheta 
        \Big\|
                \tfrac{|X(s) - Y_s|} 
                {\exp (\int_0^s \frac12\chi_r {\rm d}r)}
        \Big\|^2_{L^{p}(\Omega; \mathbb{R})}
        {\rm d}s
    +
        \tfrac{1}{8}
            \sup_{t\in [0,\vartheta]} 
            \Big\|
            \tfrac{|X(t)-Y_t|}{\exp (\int_0^t \frac{1}{2}{\chi_r} {\rm d}r)}
            \Big\|^2_{L^p(\Omega;\mathbb{R})}
            +Ch^2,
        \end{aligned}
      \end{equation}
which together with
\eqref{eq:order_one_estima_origin}
in turns gives, for any $\vartheta\in[0,T]$,
\begin{equation}
\begin{aligned}
    &
    \Big\| \sup_{t\in[0,\vartheta]}
            \tfrac{|X(t) - Y_t|^2}
            {\exp (\int_0^t \chi_s {\rm d}s)}
    \Big\|_{L^{\frac{p}{2}}(\Omega; \mathbb{R})}
    =
    \Big\| \sup_{t\in[0,\vartheta]}
            \tfrac{|X(t) - Y_t|}
            {\exp (\int_0^t \frac12 \chi_s {\rm d}s)}
    \Big\|^2_{L^{p}(\Omega; \mathbb{R})}
\\
&\leq 
C
    \int_0^\vartheta
        \Big\|
        \sup_{s\in[0,t]}
                \tfrac{|X(s) - Y_s|} 
                {\exp (\int_0^s \frac12\chi_r {\rm d}r)}
        \Big\|^2_{L^{p}(\Omega; \mathbb{R})}
        {\rm d}t
    +
        \tfrac{1}{8} 
            \Big\|
            \sup_{t\in [0,\vartheta]}
            \tfrac{|X_t-Y_t|}{\exp (\int_0^t \frac{1}{2}{\chi_r} {\rm d}r)}
            \Big\|^2_{L^p(\Omega;\mathbb{R})}
            +Ch^2.
        \end{aligned}
      \end{equation}
So the Gronwall inequality promises
\begin{equation}
\begin{aligned}
    &\Big\| \sup_{t\in[0,T]}
            \tfrac{|X(t) - Y_t|^2}
            {\exp (\int_0^t \chi_s {\rm d}s)}
    \Big\|_{L^{\frac{p}{2}}(\Omega; \mathbb{R})} 
\leq 
Ch^2 e^{CT}.
        \end{aligned}
      \end{equation}
Finally, we rely on H\"{o}lder's inequality to derive
\begin{equation}
    \begin{aligned}
            \!
            \Big\|
                \sup_{t\in[0,T]} |X(t)-Y_t|
                \Big\|_{L^{\theta}(\Omega;\mathbb{R})}
        \leq    
            \Big\|\sup_{t\in[0,T]}
                \tfrac{|X(t) - Y_t|}
                {\exp (\int_0^t {\frac{1}{2}\chi_r} {\rm d}r)} 
            \Big\|_{L^{p}(\Omega;\mathbb{R})}
            \! \cdot 
            \Big\| {\exp \Big(\int_0^{T} {\tfrac{1}{2}{\chi_r}} {\rm d}r
            \Big)}\Big\|_{L^q(\Omega;\mathbb{R})}
        \leq C_Th,
        \end{aligned}
    \end{equation}
and hence finish the proof.\qed

\begin{remark} \label{rem:strong-weak}
As a by-product of Theorem \ref{thm:optimal_converge_rate}, one can obtain the order one weak convergence of the SSAV method in finite time horizon, but with the error constant exponentially depending on the time length $T$. Such weak convergence results also hold for test functions with exponential growth, due to the  exponential integrability properties \eqref{eq:expon_integ_prope_exact_solution}  and \eqref{eq:expon_integ_prope_numerical_solution}, which may be of particular interest in practical applications.
\end{remark}
The next proposition reveals the error between the exact energy and the numerical energy.
\begin{proposition}\label{propo:energy_error}
    Under the assumptions of Theorem \ref{thm:optimal_converge_rate}, we have
\begin{equation}
    \sup_{n\in\{0,1,...,N_T\}}
    \mathbb{E}
    \Big[
    \big|
    \mathcal{H}(v(t_n),u(t_n),\rho(t_n))
    -
    \mathcal{H}(v_n,u_n,\rho_n)
    \big|^p
    \Big]
    \leq
    C_Th^p,
    \quad
    h\in(0,1],
\end{equation}
where $C_T$ is some constant independent of $h$ (exponentially depends on $T$).
\end{proposition}
\textbf{Proof:}
By virtue of Theorem \ref{thm:optimal_converge_rate} we derive
\begin{equation}
\begin{aligned}
&\mathbb{E}
    \big[
    \big|
    \mathcal{H}(v(t_n),u(t_n),\rho(t_n))
    -
    \mathcal{H}(v_n,u_n,\rho_n)
    \big|^p
    \big]
\\
&\quad =
\mathbb{E}
    \big[
    \big|
   \tfrac{1}{2}\langle v(t_n)+v_n, v(t_n)-v_n \rangle
   +
   \alpha 
   \langle u(t_n)+u_n, u(t_n)-u_n \rangle
   +
   (\rho(t_n)+\rho_n)(\rho(t_n)-\rho_n) 
    \big|^p
    \big]
\\
&\quad \leq
C 
\Big(
\mathbb{E}
    \big[
    |v(t_n)+v_n|^p |v(t_n)-v_n|^p
    +
    |u(t_n)+u_n|^p |u(t_n)-u_n|^p
    \big]
    +
    |\rho(t_n)+\rho_n|^p |\rho(t_n)-\rho_n|^p
    \big]
\Big)
\\
&\quad \leq 
C 
\big(
\mathbb{E}
    \big[
    |v(t_n)-v_n|^{2p}
    \big]
\big)^{\frac{1}{2}}
    +
C 
\big(
\mathbb{E}
    \big[
    |u(t_n)-u_n|^{2p}
    \big]
\big)^{\frac{1}{2}}
    +
C 
\big(
\mathbb{E}
    \big[
    |\rho(t_n)-\rho_n|^{2p}
    \big]
\big)^{\frac{1}{2}}
\\
&\quad \leq Ch^{p}
+
C 
\big(
\mathbb{E}
    \big[
    |\rho(t_n)-\rho_n|^{2p}
    \big]
\big)^{\frac{1}{2}}
\end{aligned}
\end{equation}
for any $n \in \{0,1,...,N_T\}$,
where one observes
\begin{equation}\label{eq:error_rho_t_n_to_rho_n}
\begin{aligned}
        \rho(t_n)-\rho_n
        &=
        \sqrt{\kappa \Phi(u(s))+C_H-\alpha |u(s)|^2}
        -
        \sqrt{\kappa \Phi(u_n)+C_H-\alpha |u_n|^2}
\\
&\quad +\sqrt{\kappa \Phi(u_n)+C_H-\alpha |u_n|^2}-\rho_n.
\end{aligned}
\end{equation}
On the one hand, due to Theorem \ref{thm:optimal_converge_rate},
it is easy to check
\begin{equation}
    \mathbb{E}
    \Big[
    \big|
    \sqrt{\kappa \Phi(u(s))+C_H-\alpha |u(s)|^2}
        -
        \sqrt{\kappa \Phi(u_n)+C_H-\alpha |u_n|^2}
    \big|^p
    \Big]
    \leq Ch^p.
\end{equation}
On the other hand, recalling the claim \eqref{eq:claim_rho_k_Order_h} gives
\begin{equation}
\begin{aligned}
&
\mathbb{E}
\Big[
\big|
\sqrt{\kappa \Phi(u_n)+C_H-\alpha |u_n|^2}
    -\rho_n
\big|^p
\Big]
\leq
Ch^p.
\end{aligned}
\end{equation}
The desired result follows, by taking the above estimates into account.
\qed

\section{Approximations of invariant measures}\label{section:long_time_weak}
\subsection{Settings}
In this section,
we pursue the analysis for
long-time behaviors of the SSAV method \eqref{eq:SSAV_unsolve_equation}-\eqref{eq:SSAV_OU_equation}. To begin with, for $x=(v,u)\in \mathbb{R}^{2m}$, we define
\begin{equation}
    \Theta(v,u):=
    \tfrac{|v|^2}{2}+\kappa \Phi(u)+
    \tfrac{\gamma}{2}\langle v,u \rangle
    +\tfrac{\gamma^2}{4}|u|^2+1
\end{equation}
and for $X^x(t),t \geq 0$ being the solution of \eqref{eq:langevin_sde} that initiates at $x\in \mathbb{R}^{2m}$,
\begin{equation}\label{eq:def_Lambda_t_x}
    \Lambda(t,x):=\mathbb{E}\big[\varphi(X^x(t))\big],
    \
    t\geq 0,
\end{equation}
where $\varphi \in C^2(\mathbb{R}^{2m}; \mathbb{R})$ with its first and second order derivatives admitting polynomial growth.
Under Assumption \ref{ass:strong_solution},
$\Lambda(t,x)$ is 
twice differentiable spatially
such that for any $\xi_1 \in \mathbb{R}^m$, 
\begin{equation}
    D\Lambda(t,x)\xi_1=\mathbb{E}\big[D\big( \varphi(X^x(t))\big)\big]
    =
    \mathbb{E}\big[\big\langle D\varphi(X^x(t)), DX^x(t)\xi_1
    \big\rangle
    \big]
\end{equation}
and for any $\xi_1,\xi_2 \in \mathbb{R}^m$, 
\begin{equation}
    D^2\Lambda(t,x)(\xi_1,\xi_2)
    =
    \mathbb{E}\big[D^2\varphi(X^x(t))
    \big( DX^x(t)\xi_1, DX^x(t)\xi_2\big)
    +
    \big\langle 
    D \varphi(X^x(t)),
     D^2 X^x(t)(\xi_1,\xi_2)
    \big\rangle
    \big].
\end{equation}
Moreover, $\Lambda(t,x)$ solves the following parabolic problem (see, e.g., \cite[Theorem 1.6.2]{cerrai2001second} and \cite{brehier2023approximation}),
\begin{equation}\label{eq:kolmo_equation}
    \tfrac{\partial \Lambda}{\partial t}(t,x)
    =\langle  D\Lambda(t,x), F(x)\rangle 
    +\tfrac{1}{2}
    \operatorname{tr}\big(GG^*D^2\Lambda(t,x)\big),
    \
    t>0;
    \
    \Lambda(0,x)=\varphi(x).
    \end{equation}
\subsection{Weak convergence and approximation of invariant measure}
As indicated by Lemma \ref{lem:energy_law} and Corollary \ref{corol:num_energy_law}, the evolution of total energy for both the exact and numerical solutions implies that their second order moments grow at most linearly in time $t$, which is evidently better than those obtained from  Corollary \ref{coro:p_order_and_expone_moment}, which grow exponentially. Moreover, we can obtain moment bounds of higher order in the following proposition.
\begin{proposition}\label{prop:bounded_moment_exact_numer_solu}
    Let Assumptions \ref{ass:strong_solution}, \ref{ass:well_pose_rho_t} hold and let $\{X(t) ,
    t\geq 0\}$ and the SSAV scheme $\{Y_{t_n}, n\in \mathbb{N}_0\}$ be defined by \eqref{eq:define_X_t_and_Y_n}.
    For any $p\in \mathbb{N}$ and $T>0$, there exists a constant $C_p$ independent of $T, h$ such that
    \begin{equation}
        \begin{aligned}
      \sup_{t \in [0,T]}
      \mathbb{E}[|X(t)|^p] {\ \textstyle \bigvee}
      \sup_{h\in(0,1]}
       \sup_{n\in\{0,1,...,N_T\}}
        \mathbb{E}[|{Y}_{t_n}|^p+|{\rho}_n|^p]  
        \leq C_p\big(\mathbb{E}[|X_0|^p]+T^{\frac{p}{2}}\big).
      \end{aligned}
    \end{equation}
\end{proposition}
\textbf{Proof:}
In view of Lemma \ref{lem:energy_law},
we first note that
    \begin{equation}
    \begin{aligned}
    H(v(t),u(t))-H(v_0,u_0)
    &=
    \int_0^t
    \big(
    -\gamma |v(s)|^2 + \tfrac{1}{2}\|\Gamma\|^2
    \big)
    {\rm d}s
    +
    \int_0^t  
        v(s)^*\Gamma 
{\rm d}{W}_s. 
    \end{aligned}
    \end{equation} 
Since $H(v(t),u(t))>0$ for all $t >0$, 
the It\^{o} formula then gives for any $p \geq 1$,
\begin{equation}\label{eq:high_ord_poly_moments_ito}
\begin{aligned}
    &H^p(v(t),u(t)) - H^p(v_0,u_0)
\\
    &\quad = \int_0^t 
    \Big(
    pH^{p-1} (v(s),u(s))     
    \big(
    -\gamma |v(s)|^2 + \tfrac{1}{2}\|\Gamma\|^2 
    \big)
    +\tfrac{p(p-1)}{2}H^{p-2}(v(s),u(s)) 
    |v(s)^*\Gamma|^2
    \Big)
     {\rm d}s
\\
&\qquad 
+
\int_0^t 
    pH^{p-1} (v(s),u(s))     
    v(s)^*\Gamma
     {\rm d}W_s
\\
&\quad \leq 
\int_0^t 
    \Big(
    \tfrac{p}{2}H^{p-1} (v(s),u(s))     
     \|\Gamma\|^2 
    +\tfrac{p(p-1)}{2}H^{p-2}(v(s),u(s)) 
    |v(s)|^2\|\Gamma\|^2
    \Big)
     {\rm d}s
\\
&\qquad 
+
\int_0^t 
    pH^{p-1} (v(s),u(s))     
    v(s)^*\Gamma
     {\rm d}W_s,
\end{aligned}
\end{equation}
which,
by recalling the definition \eqref{eq:hamiltonian_quantity} and taking expectation on both sides, implies
\begin{equation}
\begin{aligned}
    &\mathbb{E}\big[H^p(v(t),u(t))\big] -\mathbb{E}\big[ H^p(v_0,u_0)\big]
\leq 
\int_0^t 
    {p(p-\tfrac{1}{2})}\|\Gamma\|^2 \mathbb{E}\big[ H^{p-1}(v(s),u(s))  \big]
     {\rm d}s.
\end{aligned}
\end{equation}
Taking $p=2$ and using Lemma \ref{lem:energy_law}, we arrive at 
\begin{equation}
\begin{aligned}
    \mathbb{E}\big[H^2(v(t),u(t))\big] -\mathbb{E}\big[ H^2(v_0,u_0)\big]
&\leq 
\int_0^t 
    3\|\Gamma\|^2 \mathbb{E}\big[ H(v(s),u(s))  \big]
     {\rm d}s
\\
&\leq 
\int_0^t 
3\|\Gamma\|^2
\big(
\mathbb{E}\big[H(v_0,u_0)\big]+\tfrac{1}{2}\|\Gamma\|^2s
\big)
{\rm d}s
\\
&=
3\|\Gamma\|^2 \mathbb{E}\big[H(v_0,u_0)\big] t
+\tfrac{3}{4}\|\Gamma\|^4t^2.
\end{aligned}
\end{equation}
Evidently, there exists a constant $C_2$ such that
\begin{equation}
\begin{aligned}
    &\mathbb{E}\big[H^2(v(t),u(t))\big] \leq C_2\big(\mathbb{E}\big[ H^2(v_0,u_0)\big] +t^2\big).
\end{aligned}
\end{equation}
By iteration one has
\begin{equation}
\begin{aligned}
    &\mathbb{E}\big[H^p(v(t),u(t))\big] \leq C_p\big(\mathbb{E}\big[ H^p(v_0,u_0)\big] +t^p\big),
\end{aligned}
\end{equation}
which together with the fact
$$
|X(t)|^p \leq C H^{\tfrac{p}{2}}(v(t),u(t))
$$
gives the required boundedness of the exact solution.
Likewise, recalling \eqref{eq:def_v_s_n} and concerning the SSAV scheme \eqref{eq:SSAV_unsolve_equation}-\eqref{eq:SSAV_OU_equation}, we have, for any $t\in[t_n,t_{n+1}]$,
\begin{equation}
\begin{aligned}
\mathcal{H}({v}_{t,n}, {u}_{n+1}, {\rho}_{n+1})
&=
\mathcal{H}({v}_{t_n,n}, {u}_{n+1}, {\rho}_{n+1})
+
 \int_{t_n}^{t}
            \big(
            -\gamma |{v}_{s,n}|^2
            +\tfrac{1}{2}\|\Gamma\|^2
            \big)
        {\rm d}s
        +
        \int_{t_n}^{t}
            {v}_{s,n}^* \Gamma
        {\rm d}W_s. 
\end{aligned}
\end{equation}
In a similar manner as \eqref{eq:high_ord_poly_moments_ito} and by Theorem \ref{thm: numer_energy_preserv}, we can show for $p\geq 1$,
\begin{equation}
\begin{aligned}
&
\mathbb{E}\big[\mathcal{H}^p({v}_{t,n}, {u}_{n+1}, {\rho}_{n+1})\big]
-
\mathbb{E}\big[\mathcal{H}^p({v}_n, {u}_{n}, {\rho}_{n})\big]
\\
&\quad 
=\mathbb{E}\big[\mathcal{H}^p({v}_{t,n}, {u}_{n+1}, {\rho}_{n+1})\big]
-
\mathbb{E}\big[\mathcal{H}^p({v}_{t_n,n}, {u}_{n+1}, {\rho}_{n+1})\big]
\\
&\quad \leq
 \int_{t_n}^{t}
    p(p-\tfrac{1}{2})\|\Gamma\|^2
    \mathbb{E}\big[\mathcal{H}^{p-1}({v}_{s,n}, {u}_{n+1}, {\rho}_{n+1})\big]
        {\rm d}s.
\end{aligned}
\end{equation}
Taking $p=2$ and proceeding in the same way as in
\eqref{eq:evolution_energy_proof_SSAVA1},
one deduces
for $s\in[t_n,t_{n+1}]$ that
\begin{equation}
\begin{aligned}
    \mathbb{E}\big[\mathcal{H}({v}_{s,n}, {u}_{n+1}, {\rho}_{n+1})\big]
    \leq  
    \mathbb{E}\big[\mathcal{H}({v}_{n}, {u}_{n}, {\rho}_{n})\big]+\tfrac{1}{2}\|\Gamma\|^2h
\leq 
    \mathbb{E}\big[
    \mathcal{H}(v_0,u_0,\rho_0)
    \big]
    +
    \tfrac{1}{2}\|\Gamma\|^2t_{n+1}.
\end{aligned}
\end{equation}
Hence we get
\begin{equation}
\begin{aligned}
&
\mathbb{E}\big[\mathcal{H}^2({v}_{n+1}, {u}_{n+1}, {\rho}_{n+1})\big]
\\
&\quad \leq 
\mathbb{E}\big[\mathcal{H}^2({v}_n, {u}_{n}, {\rho}_{n})\big]
+
 p(p-\tfrac{1}{2})\|\Gamma\|^2 h
 \Big(
 \mathbb{E}\big[
    \mathcal{H}(v_0,u_0,\rho_0)
    \big]
    +
    \tfrac{1}{2}\|\Gamma\|^2t_{n+1}
\Big)
\\
&\quad \leq 
\mathbb{E}\big[\mathcal{H}^2(v_0, u_0, \rho_{0})\big]
+
 p(p-\tfrac{1}{2})\|\Gamma\|^2
 t_{n+1}
 \Big(
 \mathbb{E}\big[
    \mathcal{H}(v_0,u_0,\rho_0)
    \big]
    +
    \tfrac{1}{2}\|\Gamma\|^2t_{n+1}
\Big)
\\
&\quad \leq 
C_2
\Big(
\mathbb{E}\big[\mathcal{H}^2(v_0, u_0, \rho_{0})\big]
+
 t^2_{n+1}
\Big),
\end{aligned}
\end{equation}
which proves the moment bounds of the SSAV scheme for $p=2$.
Repeating the process above one
concludes the results for $p > 2$ and thus finishes the proof.
\qed

To analyze the long-time behaviors between the exact solution and the numerical discretization, some dissipative-type assumptions are necessary. However, no convexity assumptions are imposed throughout this paper.
The next assumption can be  found in \cite[Condition 3.1]{mattingly2002ergodicity}.

\begin{assumption}\label{ass:dissipative_assum_Phi}
    Assume $\Phi(x) \in C^{\infty}(\mathbb{R}^m;\mathbb{R})$ and there exist  $\beta_1>0,\beta_2\in(0,1)$ such that
    \begin{equation}
        \tfrac{1}{2}\langle \nabla \Phi(x),x \rangle
        \geq 
        \beta_2 \Phi(x)
        +\gamma^2 \tfrac{\beta_2(2-\beta_2)}{8(1-\beta_2)}|x|^2-\beta_1.
    \end{equation}
\end{assumption}
Assumption \ref{ass:dissipative_assum_Phi} guarantees the exponential ergodicity of \eqref{eq:langevin_sde}, as stated below (cf. \cite[Theorem 3.2]{mattingly2002ergodicity}).
\begin{lemma}\label{lem:expon_conver_invar_measure}
    Under Assumptions \ref{ass:well_pose_rho_t}, \ref{ass:dissipative_assum_Phi}, the exact solution $\{X(t)\}_{t\geq0 }$ of \eqref{eq:langevin_sde} admits a unique invariant measure $\mu_{\infty}$ on $\mathbb{R}^{2m}$. Moreover, for $x_0\in \mathbb{R}^{2m}$ and any measurable function $\varphi \colon \mathbb{R}^{2m} \rightarrow \mathbb{R}$ satisfying $|\varphi(x)|\leq \Theta(x)^{\ell}, \ell\geq 1$, there exist some constants $C_\ell, \lambda_\ell >0$ such that 
    \begin{equation}
        \Big|
        \mathbb{E}[\varphi(X^{x_0}(t))]- \int_{\mathbb{R}^{2m}} \varphi {\rm d}\mu_{\infty}
        \Big|
        \leq C_\ell\Theta(x_0)^\ell e^{-\lambda_\ell t},\ \forall t\geq 0.
    \end{equation}
\end{lemma}
The next assumption, combined with Assumption \ref{ass:well_pose_rho_t}, ensures the essential exponential decay property (Lemma \ref{lem:expon_decay_property} below) of the kinetic Langevin
dynamics \eqref{eq:langevin_sde}; see \cite[Hypothesis 1.1, Theorem 3.1]{talay2002stochastic} for more details.
\begin{assumption}
\label{ass:expone_decay}
Assume $\Phi(x) \in C^{\infty}(\mathbb{R}^m;\mathbb{R})$ with $\Phi(\cdot)$ and all its derivatives being of polynomial growth. Suppose that there exist some constants $0<\vartheta, 0<r\leq \gamma \leq \tilde{C}$ and function $R:\mathbb{R}^{2m} \rightarrow \mathbb{R}$ with second derivatives having polynomial growth such that
\begin{itemize}
    \item 
    $
    \vartheta(|v|^{r}+|u|^{r})
    \leq
    H(v,u)+R(v,u)+\widetilde{C},
    $
    \item
    $\mathcal{L} H(v,u)+
     \mathcal{L} R(v,u)
     \leq -\vartheta\big(H(v,u)+R(v,u)\big)+\widetilde{C},
    $
    \item 
    $
    |\partial_v R(v,u)|^2 \leq
    \widetilde{C}
    \big(
    1+H(v,u)+R(v,u)
    \big),
    $
\end{itemize}
where $v,u\in \mathbb{R}^m$ and
$\mathcal{L}\psi(v,u):=\langle v, \partial_u \psi \rangle -
\langle 
\kappa \nabla \Phi(u)+\gamma v
, 
\partial_v \psi \rangle
+\tfrac{1}{2}\sum_{i=1}^m
\partial_{v_iv_i} \psi
$
for $\psi \in C^2(\mathbb{R}^{2m};\mathbb{R})$.
\end{assumption}
It is not difficult to check that the Gaussian mixture potential \eqref{eq:Gaussian_mixture_potential}  and the double-well potential \eqref{eq:Double-well_potential}
satisfy Assumptions \ref{ass:dissipative_assum_Phi} and \ref{ass:expone_decay}.
\begin{lemma}\label{lem:expon_decay_property}
    Let Assumption \ref{ass:well_pose_rho_t} and Assumption \ref{ass:expone_decay} hold. Then there exist some constants $r_1,r_2\geq 0$ and $\eta_1,\eta_2,C>0$, possibly depending on $\varphi$ but independent of $t$, such that
    \begin{equation}
        |D\Lambda(t,x)\xi_1|
        \leq C(1+|x|^{r_1})e^{-\eta_{1} t}|\xi_1|,
        \quad
        t>0, \ \xi_1 \in\mathbb{R}^m
    \end{equation}
and 
\begin{equation}
    |D^2\Lambda(t,x)(\xi_1,\xi_2)|
    \leq
     C(1+|x|^{r_2})e^{-\eta_{2} t}|\xi_1||\xi_2|, \quad t>0, \ \xi_1,\xi_2\in\mathbb{R}^m.
    \end{equation}
\end{lemma}

We are now ready to present the main result of this section.

\begin{theorem}[Convergence for approximation of invariant measure]
\label{thm:longtime_weak_conver}
    Let Assumptions \ref{ass:well_pose_rho_t}, \ref{ass:dissipative_assum_Phi} and \ref{ass:expone_decay} hold. There exist some constants  $C_1,C_2,\lambda> 0, l\geq 0$ such that for any $t_N = Nh >0$, $N \in \mathbb{N}$, the SSAV method \eqref{eq:SSAV_unsolve_equation}-\eqref{eq:SSAV_OU_equation} admits 
\begin{equation}
\label{eq:main-result-measure-approx}
    \begin{aligned}
    &\Big|
            \mathbb{E}\big[
            \varphi(Y^{x_0}_{t_N})
                      \big]
            -
            \int_{\mathbb{R}^{2m}} \varphi {\rm d}\mu_{\infty}
    \Big|
&\leq
C_1
e^{-\lambda t_N}
+
C_2(1+(t_N)^{l})h.
    \end{aligned}
\end{equation} 
\end{theorem}
\textbf{Proof:} Using the It{\^o} formula and recalling \eqref{eq:kolmo_equation}, we acquire
\begin{equation}
    \begin{aligned}
        &\big|
            \mathbb{E}\big[
            \varphi(X^{x_0}(t_N))
                      \big]
            -
            \mathbb{E}\big[
                \varphi(Y^{x_0}_{t_N})
                    \big]
         \big|
         =
            \big|
            \mathbb{E}\big[\Lambda(t_N,x_0)\big]
            -\mathbb{E}\big[\Lambda(0,Y^{x_0}_{t_N})\big]
            \big|
        \\
                &\quad = \Big|\sum\limits_{n=0}^{N-1}\mathbb{E}\big[\Lambda(t_N-t_{n+1},Y^{x_0}_{t_{n+1}})-\Lambda(t_N-t_n,Y^{x_0}_{t_n})\big]\Big|
        \\
                &\quad = \Big|\sum\limits_{n=0}^{N-1}\mathbb{E}\Big[
                \int_{t_n}^{t_{n+1}} 
                -\tfrac{\partial \Lambda}{\partial t}(t_N-s,Y^{x_0}_s)
                    +D\Lambda(t_N-s,Y^{x_0}_s)a(s)
                    +\tfrac{1}{2} \operatorname{tr}\big(bb^*
                    D^2\Lambda(t_N-s,Y^{x_0}_s)\big) 
                    {\rm d}s
                \Big]\Big|
        \\
                &\quad = \Big|\mathbb{E}\Big[
                \int_{0}^{t_{N}} 
                D\Lambda(t_N-s,Y^{x_0}_s)a(s)
                -D\Lambda(t_N-s,Y^{x_0}_s)F(Y^{x_0}_s)
        \\
                & \qquad\qquad\quad 
                +\tfrac{1}{2} \operatorname{tr}\big(bb^*
                    D^2\Lambda(t_N-s,Y^{x_0}_s)\big)
                    -
                    \tfrac{1}{2} \operatorname{tr}\big(GG^*
                    D^2\Lambda(t_N-s,Y^{x_0}_s)\big)
                {\rm d}s
                \Big]\Big|
        \\
            &\quad = \Big|\mathbb{E}\Big[
                \int_{0}^{t_{N}} 
                D\Lambda(t_N-s,Y^{x_0}_s)
                \big(
                a(s)
                -F(Y^{x_0}_s)
                \big)
                {\rm d}s
                \Big]\Big|,
    \end{aligned}
\end{equation}
which, 
by \eqref{eq:decompo_f_Y_s-a_s},
can be expanded to give
\begin{equation}\label{eq:long-time_S_1_S_2_S_3}
    \begin{aligned}
    &\big|
            \mathbb{E}\big[
            \varphi(X^{x_0}(t_N))
                      \big]
            -
            \mathbb{E}\big[
                \varphi(Y^{x_0}_{t_N})
                    \big]
         \big|
    \\
&\quad \leq 
        \Big|\mathbb{E}\Big[
                \int_{0}^{t_{N}} 
                D\Lambda(t_N-s,Y^{x_0}_{\lfloor s \rfloor_h})
                \big(
                f(Y^{x_0}_s)-
                f(Y^{x_0}_{\lfloor s \rfloor_h})
                \big)
                {\rm d}s
                \Big]\Big|
\\
&\qquad +
 \Big|\mathbb{E}\Big[
                \int_{0}^{t_{N}} 
                \big(
                 D\Lambda(t_N-s,Y^{x_0}_{s})
                 -
                 D\Lambda(t_N-s,Y^{x_0}_{\lfloor s \rfloor_h})
                \big)
                \big(
                f(Y^{x_0}_s)-
                f(Y^{x_0}_{\lfloor s \rfloor_h})
                \big)
                {\rm d}s
                \Big]\Big|
\\
&\qquad
+
\Big|\mathbb{E}\Big[
                \int_{0}^{t_{N}} 
                D\Lambda(t_N-s,Y^{x_0}_{s})
                \big(
                {h^{-1}} B_{\lfloor s \rfloor_h/h,3} 
                \big)
                {\rm d}s
                \Big]\Big|
\\
&\quad=:\widetilde{S}_1+\widetilde{S}_2+\widetilde{S}_3.
    \end{aligned}
\end{equation}
To handle $\widetilde{S}_1$ properly, one recalls \eqref{eq:f(Y_s)-f(Y_n)-ItoFormula} and uses It\^{o}'s formula and Proposition \ref{prop:bounded_moment_exact_numer_solu} to infer
\begin{equation}
\begin{aligned}
    &\big|
        \mathbb{E}\big[
                D\Lambda(t_N-s,Y^{x_0}_{\lfloor s \rfloor_h})
                \big(
                f(Y^{x_0}_s)-
                f(Y^{x_0}_{\lfloor s \rfloor_h})
                \big)
                \big]
    \big|
\\
    &\quad =
    \Big|
    \mathbb{E}\Big[
                D\Lambda(t_N-s,Y^{x_0}_{\lfloor s \rfloor_h})
                \big(
                    \int_{t_n}^s 
                    \big(
                    f'(Y_r)a(r)
                    +
                    \tfrac{1}{2}\sum_{j=1}^m    
                    f''(Y_r)(be_j,be_j)
                    \big) 
                    {\rm d}r
                    \big)
                \Big]
                \Big|
\\
&\quad \leq  C e^{-\eta_1 (t_N-s)}
    \mathbb{E}
    \Big[
    (1+|Y^{x_0}_{\lfloor s \rfloor_h}|^{r_1})
    \Big|
        \int_{t_n}^s 
                    \big(
                    f'(Y_r)a(r)
                    +
                    \tfrac{1}{2}\sum_{j=1}^m    
                    f''(Y_r)(be_j,be_j)
                    \big) 
                    {\rm d}r
                \Big|
\Big]
\\
&\quad \leq 
C (1+(t_N)^{l_1}) e^{-\eta_1 (t_N-s)}h,
\end{aligned}
\end{equation}
and thus
\begin{equation}\label{eq:long-time_S_1}
\begin{aligned}
    \widetilde{S}_1
    \leq 
    C (1+(t_N)^{l_1})h
    \int_0^{t_N} e^{-\eta_1 (t_N-s)} {\rm d}s
\leq 
    C(1+(t_N)^{l_1})h.
\end{aligned}
\end{equation}
In view of Taylor's expansion and Proposition \ref{prop:bounded_moment_exact_numer_solu}, one can bound $\widetilde{S}_2$ as follows:
\begin{equation}\label{eq:long-time_S_2}
\begin{aligned}
    \widetilde{S}_2
    &= 
   \Big|\mathbb{E}\Big[
                \int_{0}^{t_{N}} 
                \int_{0}^1
                D^2\Lambda
                \big(
                t_N-s,Y^{x_0}_{{\lfloor s \rfloor_h}}+r(Y^{x_0}_{s}-Y^{x_0}_{{\lfloor s \rfloor_h}})
               \big)
               \big(
               Y^{x_0}_{s}-Y^{x_0}_{{\lfloor s \rfloor_h}},
                f(Y^{x_0}_s)-
                f(Y^{x_0}_{\lfloor s \rfloor_h})
                \big)
                {\rm d}r
                {\rm d}s
                \Big]\Big|
\\
&\leq 
C 
\int_{0}^{t_{N}} 
\int_{0}^1
e^{-\eta_2 (t_N-s)}
\mathbb{E}\Big[
\big(
1+\big|Y^{x_0}_{{\lfloor s \rfloor_h}}+r(Y^{x_0}_{s}-Y^{x_0}_{{\lfloor s \rfloor_h}})\big|^{r_2}
\big)
\big|
               Y^{x_0}_{s}-Y^{x_0}_{{\lfloor s \rfloor_h}}
\big|
\big|
                f(Y^{x_0}_s)-
                f(Y^{x_0}_{\lfloor s \rfloor_h})
\big|
\Big]
        {\rm d}r
        {\rm d}s
\\
&\leq 
 C(1+(t_N)^{l_2})h.
\end{aligned}
\end{equation}
At the moment we come to the estimate of $\widetilde{S}_3$.
Following the same lines as in estimates \eqref{eq:estimate_begin_B_K_3}-\eqref{eq:estimate_for_B_k3},
one can find some constant $\tilde{l}_3>0$ such that, for any $s\leq t_N $,
\begin{equation}
\begin{aligned}
    \big\| B_{\lfloor s \rfloor_h/h,3} \big\|_{L^{2}(\Omega;\mathbb{R}^{2m})}
    &\leq 
   C(1+(t_N)^{\tilde{l}_3}) h^2.
\end{aligned}
\end{equation}
Therefore, 
by Proposition \ref{prop:bounded_moment_exact_numer_solu}, there exists some constant $l_3>0$ such that
\begin{equation}\label{eq:long-time_S_3}
\begin{aligned}
    \widetilde{S}_3
    \leq
   C\int_{0}^{t_{N}} e^{-\eta_1 (t_N-s)}
   \mathbb{E}\big[
   \big(
   1+|Y^{x_0}_{s}|^{r_1}\big)
                \big|
                {h^{-1}} B_{\lfloor s \rfloor_h/h,3} 
                \big|  
                \big]
    {\rm d}s
\leq 
 C(1+(t_N)^{l_3})h.
\end{aligned}
\end{equation}
Gathering equations \eqref{eq:long-time_S_1_S_2_S_3}, \eqref{eq:long-time_S_1}, \eqref{eq:long-time_S_2} and 
\eqref{eq:long-time_S_3} together leads to
\begin{equation}
    \begin{aligned}
\big|
            \mathbb{E}\big[
                \varphi(Y^{x_0}_{t_N})
                    \big]
            -
            \mathbb{E}\big[
            \varphi(X^{x_0}(t_N))
                      \big]
         \big|
\leq
C(1+(t_N)^{\max\{l_1,l_2,l_3\}})h.
    \end{aligned}
\end{equation} 
This in conjunction with  Lemma \ref{lem:expon_conver_invar_measure} tells us that
\begin{equation}
    \begin{aligned}
    \Big|
            \mathbb{E}\big[
            \varphi(Y^{x_0}_{t_N})
                      \big]
            -
            \int_{\mathbb{R}^d} \varphi {\rm d}\mu_{\infty}
    \Big|
&\leq
\big|
            \mathbb{E}\big[
             \varphi(X^{x_0}(t_N))
                      \big]
            -
            \int_{\mathbb{R}^d} \varphi {\rm d}\mu_{\infty}
    \big|
+
    \big|
            \mathbb{E}\big[
                \varphi(Y^{x_0}_{t_N})
                    \big]
     -
            \mathbb{E}\big[
            \varphi(X^{x_0}(t_N))
                      \big]
         \big|
\\
&\leq
C_1
e^{-\lambda t}
+
C_2(1+(t_N)^{\max\{l_1,l_2,l_3\}})h,
    \end{aligned}
\end{equation} 
as required.
\qed

\subsection{Computational costs}
\label{subsec:computational-costs}
In this subsection, we aim to illustrate that,   weak error estimates polynomially (not exponentially) depending on the time length $t_N = Nh$ are crucial to significantly reducing the computational costs in approximations of invariant measures.
In light of Theorem \ref{thm:longtime_weak_conver}, to ensure the weak error 
%
\begin{equation} \label{eq:invariant-measure-error}
\Big|
        \mathbb{E}\big[
            \varphi(Y^{x_0}_{t_N})
                    \big]
        -
        \int_{\mathbb{R}^d} \varphi {\rm d}\mu_{\infty}
    \Big|
    < 
    \varepsilon
\end{equation}
for a given  precision $\varepsilon >0$, 
it suffices to appropriately choose $t_N$ and $h$ in \eqref{eq:main-result-measure-approx} such that
$$C_1e^{-\lambda t_N}
\leq 
\tfrac{\varepsilon}{2},
\
    C_2(1+(t_N)^{l})h
    \leq 
    \tfrac{\varepsilon}{2}.
$$
In other words, one needs to choose $t_N$ and $h$ satisfying
$$
t_N = \mathcal{O}
\big(
\log(\tfrac{1}{\varepsilon})
\big),
\
h= \mathcal{O}
\big(
\varepsilon 
\big(
\log(\tfrac{1}{\varepsilon})
\big)^{-l}
\big).
$$
The necessary number of iteration steps hence reads
\begin{equation}
\label{eq:step-number-pol}
\mathcal{N}_{pol}={t_N}h^{-1}= \mathcal{O}\big(\tfrac{1}{\varepsilon} \big( \log(\tfrac{1}{\varepsilon}) \big )^{1+l}\big).
\end{equation}
As noted in Remark \ref{rem:strong-weak}, the order-one strong convergence implies the order-one weak convergence, but with an exponential dependence on the time length $t_N$:
$$
 \big|
            \mathbb{E}\big[
                \varphi(Y^{x_0}_{t_N})
                    \big]
     -
            \mathbb{E}\big[
            \varphi(X^{x_0}(t_N))
                      \big]
         \big|
\leq C(1+e^{Lt_N})h,
$$
where $L>0$ would be a very large constant depending on the coefficients of SDEs and the test function $\varphi$.
Using such an error estimate promises
\begin{equation} 
\Big|
        \mathbb{E}\big[
            \varphi(Y^{x_0}_{t_N})
                    \big]
        -
        \int_{\mathbb{R}^d} \varphi {\rm d}\mu_{\infty}
    \Big|
    \leq C_1  e^{-\lambda t_N}
+
C ( 1 + e^{Lt_N} ) h.
\end{equation}
In this case, the necessary number of iteration steps required to achieve a given tolerance
$\varepsilon$ in \eqref{eq:invariant-measure-error} turns out to be
\begin{equation}
\label{eq:step-number-exp}
\mathcal{N}_{exp}=\mathcal{O}\big((\tfrac{1}{\varepsilon})^{1+L}\log(\tfrac{1}{\varepsilon})\big).
\end{equation}
Evidently, weak error estimates with polynomial dependence on the time length $t_N$
significantly reduce the computational
costs in approximations of invariant measures, compared to the case of exponential dependence.
Once uniform-in-time moment bounds of numerical approximations are available, which is the case for an implicit splitting scheme proposed by \cite{chen2025new}, the weak error estimate can be improved to 
\[
 \big|
            \mathbb{E}\big[
                \varphi(Y^{x_0}_{t_N})
                    \big]
     -
            \mathbb{E}\big[
            \varphi(X^{x_0}(t_N))
                      \big]
         \big|
\leq C h
\]
and the resulting error for the approximation of invariant measure in \eqref{eq:main-result-measure-approx} can hold for $l=0$. 
For this case, the necessary number of iteration steps required is thus reduced to
\begin{equation}
\label{eq:step-number-uniform}
\mathcal{N}_{uni}=\mathcal{O}\big(\tfrac{1}{\varepsilon}\log(\tfrac{1}{\varepsilon})\big).
\end{equation}
Comparing \eqref{eq:step-number-pol} with \eqref{eq:step-number-uniform}, one can clearly see, non-uniform error bounds with polynomial dependence on the time length $t_N$, instead of uniform-in-time error bounds, only result in a marginal increase of the iteration number (supplementary polynomial dependence with respect to $|\log (\varepsilon)|$). But an exponential dependence on $t_N$ would, as shown by \eqref{eq:step-number-exp}, result in a substantial increase of the cost (supplementary polynomial dependence with respect to $\varepsilon^{-1}$). 

It is worthwhile to mention that, this observation has been previously revealed in \cite{brehier2023approximation}, where an explicit tamed Euler scheme was proposed for approximation of the invariant distribution of ergodic SDEs with contractive coefficients. For the considered kinetic Langevin dynamics \eqref{eq:langevin_sde}, however, the drift fails to obey such a contractive condition.
%
Also, we point out that, only the iteration number is calculated in the above discussion, whereas the costs of every step are not even considered, which is a key component of the overall computational costs. From this point of view, the proposed explicit SSAV scheme \eqref{eq:SSAV_unsolve_equation}-\eqref{eq:SSAV_OU_equation} can be much more efficient than implicit schemes in \cite{talay2002stochastic,cui2022density,chen2025new}, particularly in the high-dimensional case $m\gg 1$.

\section{Numerical experiments} \label{section:numerical_experiment}

In this section we report some numerical experiments to confirm the theoretical findings and the effectiveness of the proposed scheme \eqref{eq:practical_iteration_SSAV} in sampling.
%
Throughout this section, we set $\alpha=1$ and $C_H=1000$ for parameters of the SSAV method.
The expectation is approximated by 
averaging over $M = 5000$ Monte Carlo sample paths. Next we consider the following three test problems.

\textbf{Example 1: Gaussian mixture potential}

Fix $m=1, \kappa =2, \gamma=1, \Gamma=2$ and for $u \in \mathbb{R}$ consider 
\begin{equation}\label{eq:Gaussian_mixture_potential}
    \Phi(u)=\tfrac{(u-\iota)^2}{2\sigma^2}
    -
    \log(\tfrac{1}{3}+
    \tfrac{2}{3}e^{-\frac{2u\iota}{\sigma^2}})
    ,\
    \nabla \Phi(u)=\tfrac{u-\iota}{\sigma^2}+\tfrac{4\iota}{\sigma^2}
    \big(
    e^{\frac{2\iota u}{\sigma^2}}+2
    \big)^{-1},
    \quad \iota \in \mathbb{R},\sigma>0
    .
\end{equation}
The corresponding invariant measure reads as: for $x = (v^*,u^*)^*\in\mathbb{R}^{2}$,
\begin{equation}\label{eq:mixture_Gaussian_measure}
    \mu_{\infty}(x)
    =
    \mu^{(1)}_{\infty}(v)
    \mu^{(2)}_{\infty}(u)
    =
    \Big(
    \tfrac{1}{ \sqrt{2\pi \kappa}}
    e^{-\frac{v^2}{2\kappa}}
    \Big)
    \Big(
    \tfrac{1}{\sqrt{2\pi\sigma^2}}
    \big(
    \tfrac{1}{3}e^{-\frac{(u-\iota)^2}{2\sigma^2}}
    +
    \tfrac{2}{3}e^{-\frac{(u+\iota)^2}{2\sigma^2}}
    \big)
    \Big).
\end{equation}
In all tests we assign $\iota=1
,\sigma=\tfrac{1}{2}$.
Clearly, $\Phi \in C^{\infty}(\mathbb{R};\mathbb{R})$ and $ \Phi'$ is Lipschitz.

\textbf{Example 2: Double-well potential}

Take $m=1, 2, 20, \kappa=1, \gamma=1,\Gamma=\sqrt{2}I_{m \times m}$ and 
\begin{equation}\label{eq:Double-well_potential}
\Phi(u)=\tfrac{1}{4}|u|^4-\tfrac{1}{2}|u|^2, u \in \mathbb{R}^m,\
\nabla \Phi(u)=|u|^2u-u,
\
u \in \mathbb{R}^m.
\end{equation}
Then its invariant measure is given by
\begin{equation}\label{eq:doublewell_example_measure}
    \mu_{\infty}(x)
    \propto
    e^{-\frac{|v|^2}{2}-(\frac{|u|^4}{4}-\frac{|u|^2}{2})},
    \
    x \in (v^*,u^*)^*.
\end{equation}

\textbf{Example 3: A bimodal distribution from \cite[Example 6.4]{rubinstein2016simulation}}

Fix $\kappa=0.1,\gamma=0.05,\Gamma=0.1$ and for $u=(u_1,u_2)^*\in\mathbb{R}^2$ we consider 
$$
\Phi(u)=\tfrac{u_1^2u_2^2+u_1^2+u_2^2-8(u_1+u_2)}{2},\ 
\nabla \Phi(u)=(u_1u_2^2+u_1-4,u_1^2u_2+u_2-4)^*.
$$

It is not difficult to check that Examples 1, 2 satisfy all the assumptions mentioned above, of which Assumption \ref{ass:expone_decay} follows from \cite[Example 1.1, Example 1.2]{talay2002stochastic} ($R_1(v,u)=\tfrac{1-c}{4}|u|^2+\tfrac{1}{2}\langle v,u\rangle$ for Example 1 and $R_2(v,u)=c\langle v,u\rangle$ for Example 2 with $c>0$ small enough). By setting $R_3(v,u)=c\langle v, u\rangle$ for some appropriate $c>0$, one can also verify that all the above assumptions are satisfied for Example 3, taken from \cite[Example 6.4]{rubinstein2016simulation}.



\begin{figure}[H]
    \centering  
        \subfigure[Strong and energy convergence rates]{
        \includegraphics[width=7cm,height = 6cm]{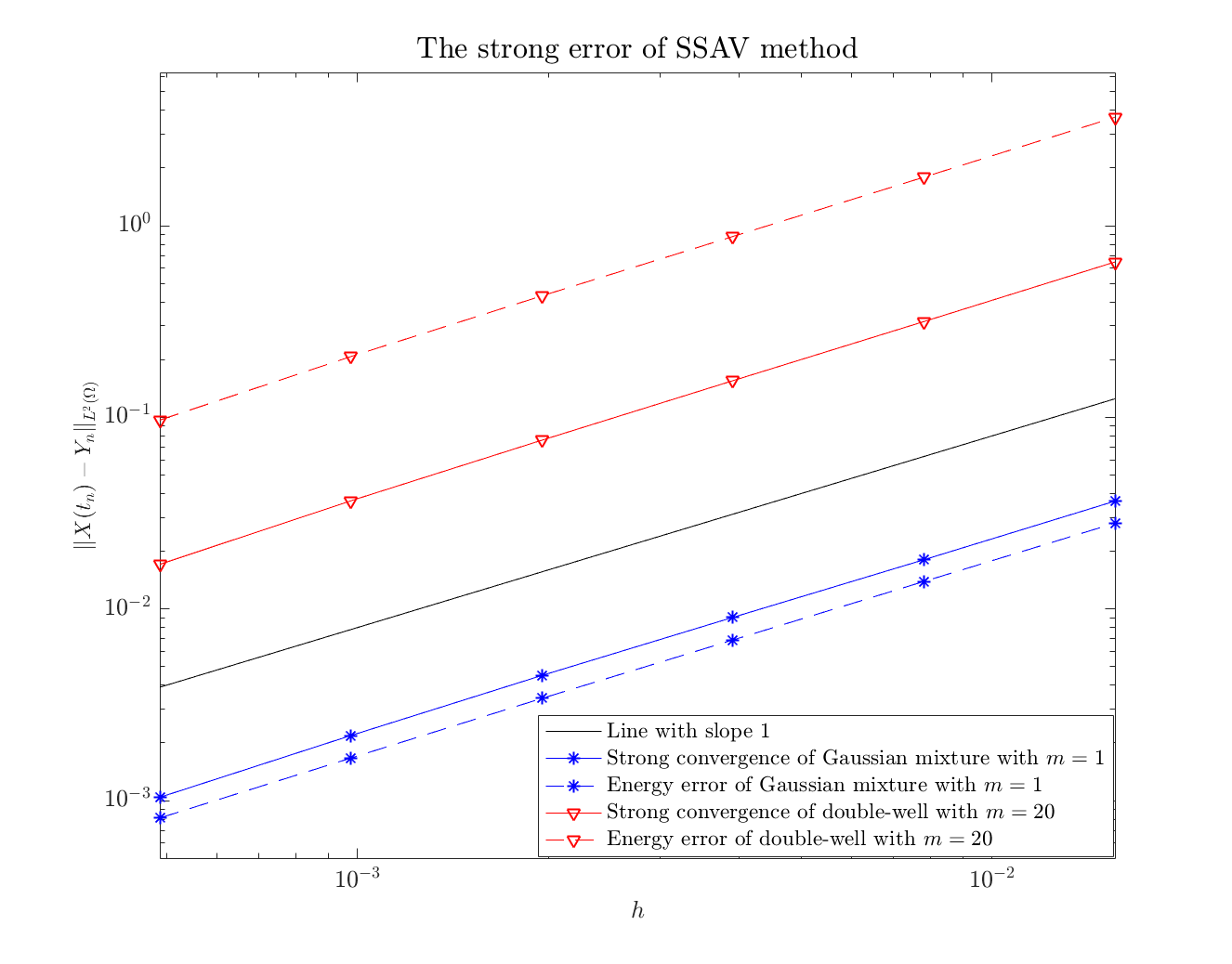}}
        \subfigure[Weak convergence rate]{
        \includegraphics[width=7cm,height = 6cm]{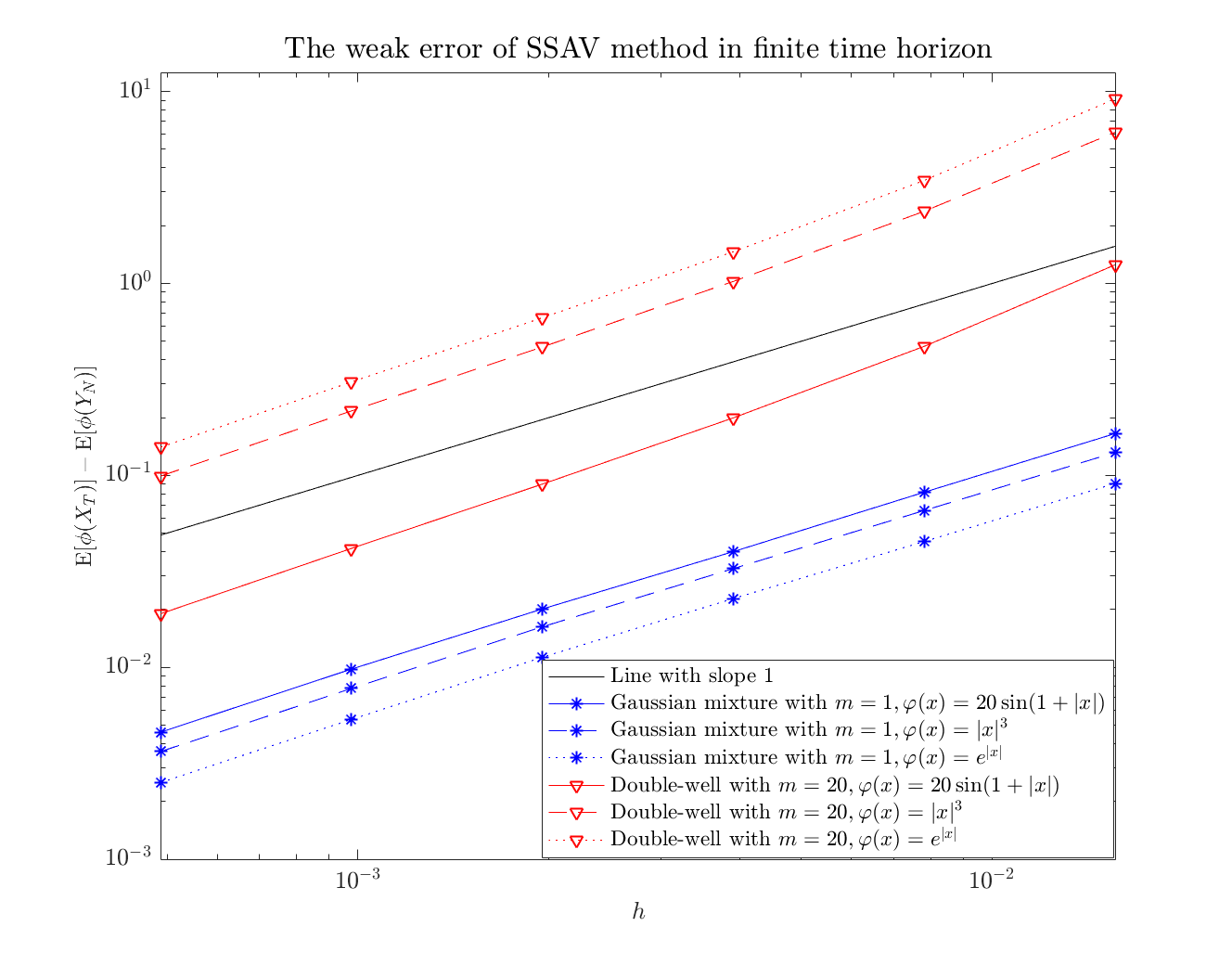}} 
        \caption{\small Strong and weak errors of the SSAV method in finite time horizon
        }
        \label{fig:strong_and_weak_error}
\end{figure}

To test strong and weak convergence rates, we let $T = 1$, $N = 2^{k}, k = 6, 7,..., 11$ and regard fine approximations with $N_{\text{exact}} = 2^{14}$ as the ``true" solution.
%
%
In terms of the mean-square error,
the strong convergence and the energy error of Example 1 and Example 2 ($m=20$) are depicted in Figure \ref{fig:strong_and_weak_error} (a) on a log-log scale,
while 
Figure \ref{fig:strong_and_weak_error} (b) displays the weak convergence 
for the test functions
$\varphi_1(x)=20\sin(1+|x|), \varphi_2(x)=|x|^3
$
and
$
\varphi_3(x)=e^{|x|}, x\in \mathbb{R}^{2m}$.
Evidently, all errors decrease at a slope close to $1$, consistent with the previous theoretical results.

\begin{figure}[H]
    \centering  
        \subfigure[Gaussian mixture potential with $m=1$]{
        \includegraphics[width=7cm,height = 6cm]{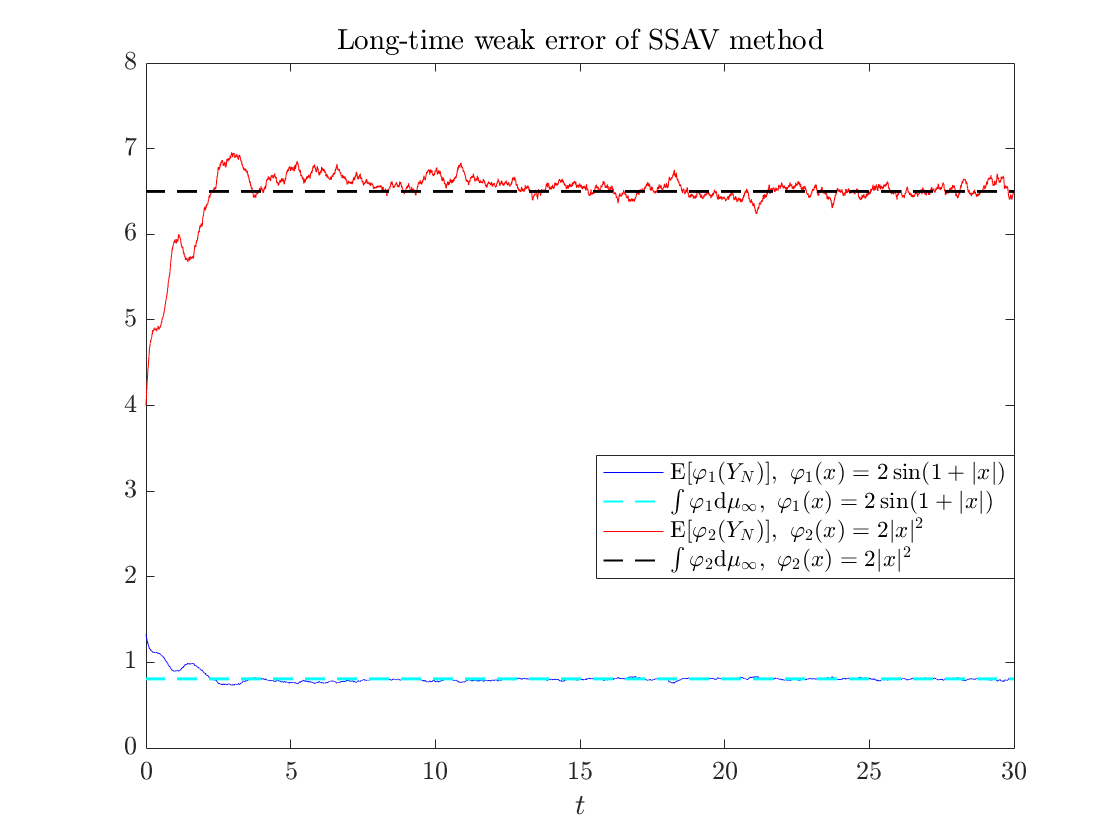}}
        \subfigure[Double-well potential with $m=1$]{
        \includegraphics[width=7cm,height = 6cm]{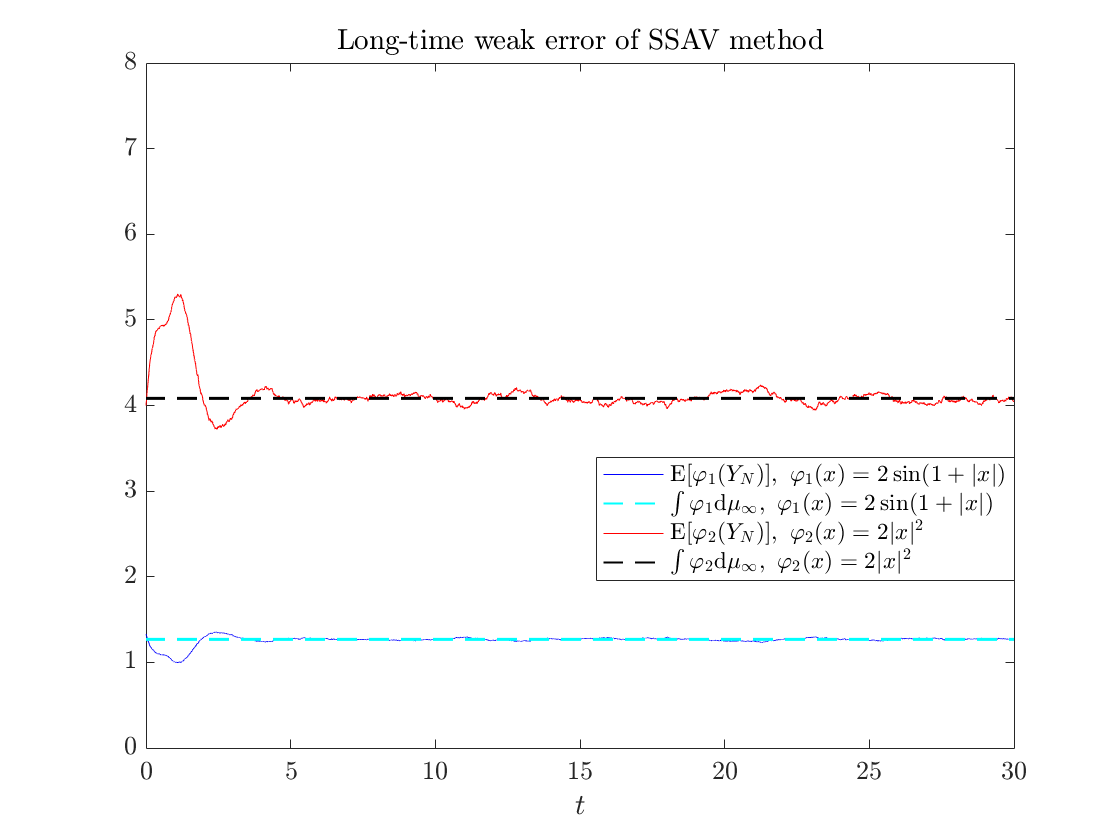}} 
        \caption{\small Long-time weak errors of the SSAV method
        }
        \label{fig:long_time_weak}
\end{figure}

Next we examine 
the long-time weak convergence of the SSAV method. Consider the  Gaussian mixture potential and double-well potential for $m=1$, along with the test functions $\varphi_1(x)=2\sin(1+|x|)$, $\varphi_2(x)=2|x|^2, x\in \mathbb{R}^2$. According to  \eqref{eq:mixture_Gaussian_measure} and \eqref{eq:doublewell_example_measure}, one can figure out the corresponding $\int_{\mathbb{R}^2} \varphi {\rm d}{\mu_{\infty}}$. Then Figure \ref{fig:long_time_weak} illustrates the evolution of error between $\mathbb{E}[\varphi(Y^{x_0}_N)]$ and $\int_{\mathbb{R}^d} \varphi {\rm d}\mu_{\infty}$ for these two models in $t\in[0,30], h=2^{-9}$.
One can observe that the SSAV method generates a good approximation after $t=5$.

\begin{figure}[H]
      \centering
      \includegraphics[scale=0.198]{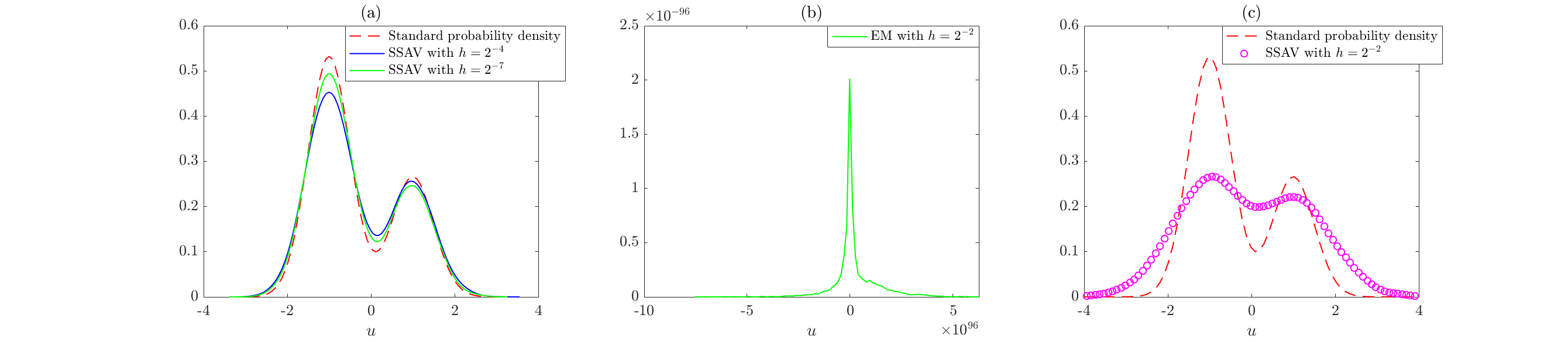}
      \caption{\small Density approximation of Gaussian mixture potential with $m=1$: (a) SSAV method with $h=2^{-4}$ and $h=2^{-7}$;  (b) EM method with $h=2^{-2}$;  (c) SSAV method with $h=2^{-2}$
      \label{fig:density_approximation_Gaussian}}
\end{figure}

\begin{figure}[H]
      \centering
      \includegraphics[scale=0.245]{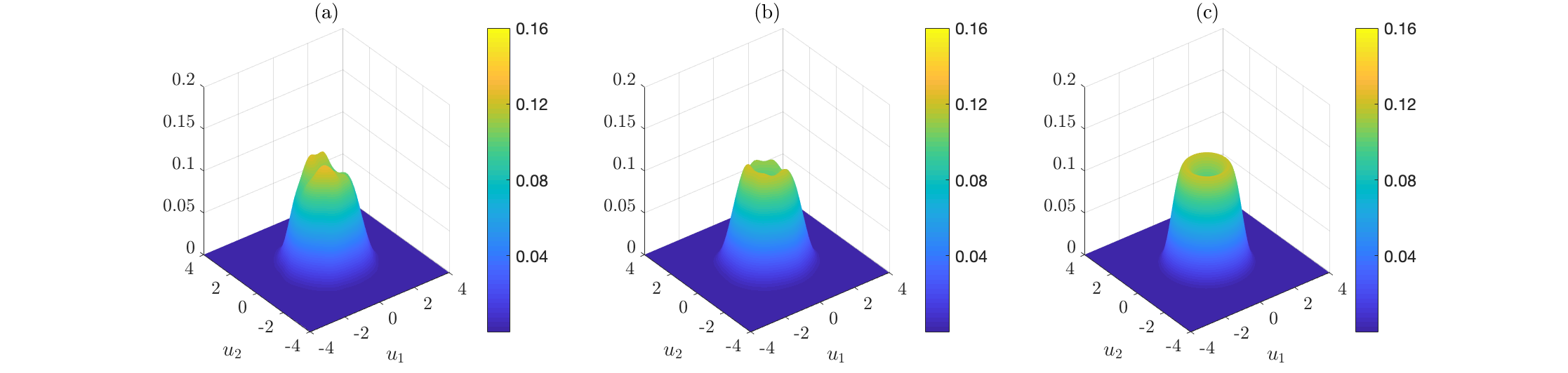}
      \caption{\small Density approximation of double-well potential with $m=2$: (a) SSAV method with $h=2^{-4}$;  (b) SSAV method with $h=2^{-7}$;  (c) Standard probability density
    \label{fig:density_approximation_doublewell}}
\end{figure}

\begin{figure}[H]
      \centering
      \includegraphics[scale=0.20]{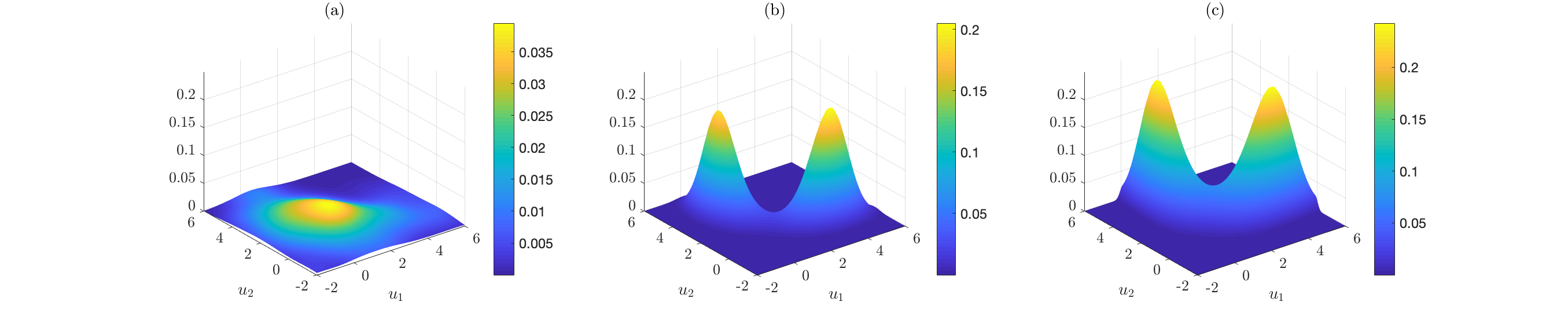}
      \caption{\small Density approximation of bimodal distribution: (a) EM method with $h=2^{-4}$; (b) SSAV method with $h=2^{-4}$;  (c) Standard probability density
    \label{fig:density_approximation_example_3}}
\end{figure}

It is well known that the kinetic Langevin
dynamics can be used for sampling from a target distribution $\pi(u)\propto \exp(-\Phi(u))$. To see the performance of SSAV method in sampling, we consider the numerical density at a large endpoint $T=500$. Figure \ref{fig:density_approximation_Gaussian} displays the density approximation for the Gaussian mixture potential with $m=1$. Specifically, Figure \ref{fig:density_approximation_Gaussian} (a) shows the density approximation using SSAV method with $h=2^{-4}$ and  $h=2^{-7}$. Figures \ref{fig:density_approximation_Gaussian} (b) and (c) compare the Euler--Maruyama (EM) method and SSAV method for sampling with the stepsize $h=2^{-2}$, where one can clearly see that the SSAV method is much more stable than the EM method when the stepsize is large. Moreover, Figure \ref{fig:density_approximation_doublewell} shows the density approximation of SSAV method for the double-well potential with $m=2$, demonstrating that our method performs well when used to sample superquadratic potentials.
Finally, we display the sampling performance of EM method and SSAV method for the bimodal distribution by using $h=2^{-4}$ and $M=20000$. As shown in Figure \ref{fig:density_approximation_example_3}, in our example, the SSAV method is significantly superior to the EM method.

\bibliographystyle{abbrv}

\bibliography{An_explicit_splitting_SAV_scheme_for_kinetic_Langevin_dynamics}
\end{document}